\documentclass[sn-mathphys]{sn-jnl}

\jyear{2021}%

\usepackage{lineno,hyperref}
\usepackage{booktabs}
\usepackage{multirow}
\usepackage{color}
\usepackage{graphicx}
\usepackage{subfigure}
\usepackage{float}
\usepackage{bm}
\usepackage{lipsum}
\usepackage{listings} 
\usepackage{amsfonts,enumerate,amsmath,amssymb}
\theoremstyle{thmstyleone}%
\def\qed{\strut\hfill $\Box$}
\newtheorem{thm}{Theorem}[section]

\newtheorem{lem}[thm]{Lemma}

\newtheorem{rem}[thm]{Remark}
\newtheorem{defn}[thm]{Definition}
\newtheorem{exm}[thm]{Example}

\newcommand{\thmref}[1]{Theorem~{\rm \ref{#1}}}
\newcommand{\lemref}[1]{Lemma~{\rm \ref{#1}}}

\def\para#1{\vskip .4\baselineskip\noindent{\bf #1}}
\numberwithin{equation}{section}
\theoremstyle{thmstyletwo}%

\begin{document}

\title[Convergence of martingale solutions to the hybrid slow-fast system]{Convergence of  martingale solutions to the hybrid slow-fast system}


\author*[1,2]{\fnm \sur{Yong Xu  }}\email{hsux3@nwpu.edu.cn}

\author[1]{\fnm \sur{Xiaoyu Yang}}\email{yangxiaoyu@yahoo.com}

\author[1]{\fnm \sur{Bin Pei}}\email{binpei@nwpu.edu.cn}

\author[3]{\fnm \sur{Yuzhen Bai}}\email{baiyu99@126.com}

\affil[1]{\orgdiv{School of Mathematics and Statistics}, \orgname{Northwestern Polytechnical University}, \orgaddress{ \city{Xi'an}, \postcode{710072}, \country{China}}}

\affil[2]{\orgdiv{MIIT Key Laboratory of Dynamics and Control of Complex Systems}, \orgname{Northwestern Polytechnical University}, \orgaddress{ \city{Xi'an}, \postcode{710072}, \country{China}}}

\affil[3]{\orgdiv{School of Mathematical Sciences}, \orgname{Qufu Normal University}, \orgaddress{\city{Qufu}, \postcode{273165}, \country{China}}}


\abstract{This paper is devoted to studying the weak convergence for a slow-fast system with jumps modulated by Markovian switching regimes with the martingale method. However, due to the coexistence of fast component and Markovian switching regimes, the martingale method and perturbed test functions can not be applied directly.  In this situation, a  combination of perturbed test functions and the time discretization is applied efficiently. And the choice of appropriate perturbed test functions, which are related to the averaged coefficients, plays a decisive role. Our results also  cover the case of slow-fast system without Markovian switching regimes. Finally, some examples are presented, and  numerical simulations are carried out to  observe a  good agreement. }




\maketitle

\section{Introduction}\label{sec1}
The multiscale models arise widely in many scientific domains, for example, suspended system \cite{Georgiou1998Slow}, mobile robot \cite{Sanfelice2011On}, aerosol particles in a cellular flow field \cite{2012Approximation} and compartmental model \cite{Krupa2008Mixed}.
In a suspended system, because the vibration frequency of springs is faster than that of tires, this leads to a multiscale system. It may still be computationally expensive for the mentioned multiscale system in the field of engineering, and it is of foremost importance to reduce the computation load. Thus, finding the low-dimensional system that can effectively describe the dynamics of the  {original} multiscale model would be an exciting practice. Fortunately, averaging principles \cite{Golec1990Averaging,Givon2006Strong} can be applied to simplify the original system by the averaged system efficiently. In other words, taking the fast-slow system, for example, it can be shown that the slow component converges to a limit averaged component in some convergence sense \cite{Xu2017Stochastic,Xu2020Averaging}. In addition, it is well known, Gaussian noise has been greatly developed to describe the stochastic disturbance in some real systems \cite{Arnold1974Stochastic}. However, Gaussian noise is not appropriate to model the stochastic disturbance in many practical problems where large jumps and bursts might exist. Fortunately, L\'evy-type processes \cite{Applebaum2009,Xu2021On}, a class of important non-Gaussian processes with jumps or bursts, are appropriate models for these special situations. 

It is worth pointing out that the aforementioned
results are focused on multiscale systems modeled by stochastic differential equations (SDEs). But, the continuous system cannot describe
the practical systems whose structures are subject to
stochastic abrupt changes caused by
abrupt phenomena \cite{Yin1998Continuous} such as component failures, repairs or other disturbances in an unreliable machine. In fact, this kind of stochastic abrupt change has been found in wireless communications, manufacturing systems and other related fields \cite{Mao2006Stochastic}. Thus, that might be a good choice to consider the continuous system modulated by Markovian switching regime (also called hybrid system) which is a dynamic system that exhibits both continuous and discrete dynamic behavior \cite{Yin2012Continuous}.  {For more details on the definitions and properties of Markovian switching, please refer \cite{Yin1998Continuous,Yin2012Continuous,Yin2010Approximation,Yin2010Hybrid}. Recently, averaging principles for stochastic partial differential equations (SPDEs) modulated by Markovian switching are investigated \cite{Pei2015Mild,Pei2018Averaging}.}

Due to the consequently increasing interest in the multiscale models with  jumps and stochastic abrupt changes
arisen from the engineering fields, it is necessary to investigate the averaging principle for slow-fast system with small jumps modulated by Markovian switching regimes as follows:
\begin{eqnarray}\label{orginal2.3}
\begin{cases}
d{x^\varepsilon}\left(t\right)&=f\left({{x^\varepsilon}\left(t\right),{y^\varepsilon}\left(t\right),{r^\varepsilon}\left({t}\right)}\right)dt+\sigma\left({{x^\varepsilon }\left(t\right),{y^\varepsilon}\left(t\right),{r^\varepsilon}\left({t}\right)} \right)d {W_1}\left(t\right)\cr
&\quad+\int_{\vert  z\vert <c}{g\left({{x^\varepsilon}\left(t\right),z}\right) {\tilde N_1}\left( {dt,dz} \right)},\cr
d{y^\varepsilon}\left(t\right)&=\frac{1}{\varepsilon } {b}\left({{y^\varepsilon}\left(t\right)}\right)dt + \frac{1}{{\sqrt \varepsilon  }} {V} \left({{y ^\varepsilon }\left( t \right)} \right)d {{W_2}}\left(t\right)\cr
&\quad +\int_{\vert z\vert <c}{ {h}\left({{y ^\varepsilon}\left(t\right),z}\right) {\tilde N_2^\varepsilon} \left({dt,dz}\right)},
\end{cases}
\end{eqnarray}
with $x^\varepsilon(0)=x_0 \in \mathbb{R}^n$ and $y^\varepsilon(0)=y_0 \in \mathbb{R}^n$, where $0 <\varepsilon \ll 1$ is a small parameter describing the ratio of the time scale between the slow component $x^\varepsilon$ and fast component $y^\varepsilon$. Let  {$\{r^\varepsilon(t)\}_{t\geq 0}$ represent two-time-scale Markovian switching regimes with a finite space $\mathcal{M}$ and $f,\sigma,g,b,V,h$ be nonlinear coefficients which will be given in Section 2.}
Let $(\Omega, \mathcal{F}, \{\mathcal{F}_t\}_{t \geq 0},\mathbb{P})$ be a stochastic basis satisfying the usual conditions,	
	 {$\{W_1\left(t\right)\}_{t \ge 0}$, $\{W_2 \left(t\right)\}_{t \ge 0}$ be independent Brownian motions} and $ {\tilde N_1\left({dt,dz} \right)}$, ${\tilde N_2^\varepsilon\left({dt,dz}\right)}$ be compensated Poisson random measures \cite{Applebaum2009}, associated with independent Poisson random measures ${N_1\left({dt,dz} \right)}$ and ${N_2\left({dt,dz}\right)}$, defined as,
$\tilde N_1\left({dt,dz}\right)=N_1\left({dt,dz}\right)-v_1\left({dz}\right)dt,$ and
$\tilde N_2^\varepsilon\left({dt,dz}\right)={N_2}\left({dt,dz}\right)-\frac{1}{\varepsilon}v_2\left({dz}\right)dt,$ respectively.

In the Brownian motion case, an approximate method for stochastic chemical kinetics with two-time
scales by chemical Langevin equations in the sense of the weak convergence was obtained in \cite{Fuke2016Approximate}. After that,  functional diffusions with two-time scales in which the fast-varying process is a rapidly-changing diffusion are examined,  then an averaging principle was obtained using weak convergence methods and perturbed test functions \cite{Wu2020An}. In the jump case,   two-time-scale jump diffusion models modulated by continuous-time Markov chains are considered and  their limit processes are derived using weak convergence methods directly \cite{Yin2004Two}.
Moreover, the weak convergence methods of Kushner \cite{Kushner1984Approximation} (see also \cite{Wu2020An,Fuke2016Approximate,Yin2004Two}) cannot be applied directly in (\ref{orginal2.3}) since we have to consider both the two-time-scale Markovian switching and jumps term. New approaches have to be developed.
Motivated by this, in this paper, we will use the techniques of perturbed test functions and time discretization to eliminate the jump term and Markovian regimes efficiently to study the convergence of martingale solution to the slow-fast system with small jumps modulated by Markovian switching regimes.  There are two main difficulties in our setting. First, due to the co-existence of fast component and two-time-scale Markovian switching regimes, the martingale method with perturbed test functions is unable to obtain the weak convergence result (see Theorem 3.1). Second, how to construct the appropriate perturbed test functions will encounter many technical difficulties.  

The rest of the paper is arranged as follows.   Section 2 presents some notations and assumptions. The main result is proved in Section 3. In Section 4, we give some examples to   demonstrate the applications of averaging. Ultimately, a conclusion is given.

\section{Slow-fast systems with two-time-scale Markov switching regimes and jumps}\label{sec-2}

In this section, some assumptions and notations are given.
$r^\varepsilon(t)$ is a fast varying Markov chains whose generator is given by 
${{\mathcal{Q}}^\varepsilon }: = \frac{{\widetilde {\mathcal{Q}}}}{\varepsilon }, $
where ${{\widetilde {\mathcal{Q}}}}/{\varepsilon }$  represents the fast-varying generator, which are bounded and Borel measurable. The generator $\widetilde {\mathcal{Q}}$ is weakly irreducible, meaning that
\begin{eqnarray*}\label{lemma3.1}
	\begin{cases}
		{\nu \widetilde  {\mathcal{Q}} = 0},\cr
		{\sum\limits_{\gamma  = 1}^n {{\nu_\gamma }}  = 1},
	\end{cases}
\end{eqnarray*}
has a unique solution $\nu = \left( {{\nu_1},{\nu_2}, \cdots ,{\nu_n}} \right)$ which is termed a quasi-stationary distribution satisfying $\nu_\gamma  \ge 0$ for each $\gamma \in \mathcal{M}=\left\{1,2, \cdots ,n\right\}$. 

Let $M$ denote the set of real-valued progressively measurable functions that are nonzero only on a bounded $t$-interval, $${M^\varepsilon }=\{ \iota\in M: \mathop {\sup }\limits_t  {\vert \iota (t) \vert}< \infty ,\iota (t): \mathcal F_t^\varepsilon{\rm -measurable}  \},$$  {where $\{{\mathcal F}_t\}_{t \ge 0}:=\sigma\{x^\varepsilon\left(s\right),y^\varepsilon (s):0\le s\le t\}$.}

Firstly we define an infinitesimal operator ${\hat {\mathcal A}^\varepsilon }$ as follows.
 {\begin{defn}\label{2-1}
		\cite{Kushner1984Approximation} We say that $\iota\left(  \cdot  \right) \in D( {{{\hat {\mathcal A}}^\varepsilon }} )$, the domain of ${\hat {\mathcal A}^\varepsilon }$, and ${\hat {\mathcal A}^\varepsilon }\iota =\kappa$, if $\kappa(t),\iota(t) \in M^\varepsilon$ and 
		\begin{eqnarray*}\label{orginal2.4}
			p\text{-}{\mathop {\lim }\limits_{\delta\to 0}\left[{\frac{{\mathbb{E}_t^\varepsilon \iota\left({t+\delta}\right)-\iota\left(t\right)}}{\delta}- \kappa\left(t\right)}\right]}=0 ,
		\end{eqnarray*}
		where  $\mathbb{E}_t^\varepsilon$ represents the  expectation conditioned on   $\mathcal F_t^\varepsilon $, and $p\text{-}\mathop {\lim }\limits_{\delta\to 0}$ is defined as: $\iota = p\text{-}\mathop {\lim }\limits_{\delta\to 0}{\iota^\delta }$ 
		\begin{eqnarray*}\label{orginal2.5}
			\begin{split}
				\begin{cases}
					\mathop {\sup }\limits_{t,\delta } \mathbb{E} \vert  {{\iota^\delta }\left( t \right)}  \vert  < \infty ,\cr
					\mathop {\lim }\limits_{\delta  \to 0} \mathbb{E} \vert  {{\iota^\delta }\left( t \right) - \iota\left( t \right)}  \vert  = 0 ,
				\end{cases} 
			\end{split}
		\end{eqnarray*}
		for each $t$. 
\end{defn}}
 {\begin{lem}\label{2-2}
		\cite{Kurtz1975Semigroups} If $\iota\left(  \cdot  \right) \in D( {{{\hat {\mathcal A}}^\varepsilon }} )$, then
		\begin{eqnarray*}\label{orginal2.6}
			\iota\left( t \right) - \int_0^t {{{\hat {\mathcal A}}^\varepsilon }\iota\left( u \right)du}  = :M_\varepsilon ^\iota\left( t \right)
		\end{eqnarray*}
		is a martingale, and also 
		\begin{eqnarray*}\label{orginal2.7}
			\mathbb{E}_t^\varepsilon \iota\left( {t + s} \right) - \iota\left( t\right) = \int_t^{t + s} {\mathbb{E}_t^\varepsilon {{\hat {\mathcal A}}^\varepsilon }\iota\left( u \right)du}
		\end{eqnarray*}
		with probability 1.
\end{lem}}
There follow additional assumptions which are needed in our argument:		
		
\begin{enumerate}[({A}1)]
	\item For each $\gamma \in \mathcal{M}$, any $ {{x_i}}  \in {\mathbb{R}^n}$ and $y_i  \in {\mathbb{R}^n}$, there exists a ${C_\gamma} > 0$, such that
	\begin{align}
	&{ \vert  {f\left( {x_1,y _1,\gamma} \right) - f\left( {x_2,y _2,\gamma} \right)}  \vert ^2}
	\vee { \vert  {\sigma \left( {x_1,y _1,\gamma} \right) - \sigma \left( {x_2,y _2,\gamma} \right)}  \vert ^2}\cr
	&\vee \int_{ \vert z  \vert  < c} {{{ \vert  {g\left( {x_1,z} \right) - g\left( {x_2,z} \right)}  \vert }^2}v_1\left( {dz} \right)} \cr
	&\le {C_\gamma}\left( {{{ \vert  {x_1 - x_2}  \vert }^2}  + {{ \vert  {y _1 - y _2}  \vert }^2}} \right).\nonumber
	\end{align}
	\item For each $\gamma \in \mathcal{M}$, all ${{x}}  \in {\mathbb{R}^n}$ and $y  \in {\mathbb{R}^n}$, i.e. there exists a positive constant ${C'_\gamma}$ such that
	\begin{align}
	\quad{ \vert  {f\left( {x,y,\gamma} \right)}  \vert ^2}
	\vee { \vert  {\sigma \left( {x,y,\gamma} \right)}  \vert ^2}
	\vee \int_{ \vert z  \vert  < c} {{{ \vert  {g\left( {x,z} \right)}  \vert }^2}v_1\left( {dz} \right)}
	\le {C'_\gamma}\left( {1+ {{ \vert  {x}  \vert }^2}  + {{ \vert  {y }  \vert }^2}} \right).\nonumber
	\end{align}
	\item There exist ${\alpha _{1}},{\alpha _{2}}, {\alpha _{3}}$, which are with $2{\alpha _{1}} - {\alpha _{2}} - {\alpha _{3}} > 0$. Then for any  ${y _1}$, ${y _2}$, we suppose that 
	\begin{align}
	\left\langle {{y _1} - {y _2},{b}\left( {{y _1}} \right) - {b}\left( {{y _2}} \right)} \right\rangle  &\le  - {\alpha _{1}}{ \vert  {{y _1} - {y _2}}  \vert ^2},\cr
	{ \vert  {{b} \left( {{y _1}} \right) - {b} \left( {{y _2}} \right)}  \vert ^2}  \vee  { \vert  {V \left( {{y _1}} \right) - V \left( {{y _2}} \right)}  \vert ^2} &\le {\alpha _{2}}{ \vert  {{y _1} - {y _2}}  \vert ^2},\cr
	\int_{ \vert  z  \vert  < c} {{{ \vert  {h\left( {{y _1},z} \right) -h\left( {{y _2},z} \right)} \vert }^2}v_2\left( {dz} \right)}  &\le {\alpha _{3}}{ \vert  {{y _1} - {y _2}}  \vert ^2}.
	\end{align}
	There exist ${\alpha '_{1}},{\alpha '_{2}},{\alpha '_{3}} > 0$   with $2{\alpha '_{1}} - {\alpha '_{2}} - {\alpha '_{3}} > 0$, such that
	\begin{align}
	\left\langle {y ,{b}\left( {y } \right)} \right\rangle &\le  - {\alpha '_{1}}{ \vert  y   \vert ^2}+\alpha ,\cr
	{ \vert  {{b} \left( {y } \right)}  \vert ^2} \vee { \vert  {V \left( {y } \right)}  \vert ^2} &\le {\alpha '_{2}}{ \vert  y   \vert ^2}+\alpha,\cr
	\int_{ \vert  z  \vert  < c} {{{ \vert  {h\left( {y ,z} \right)}  \vert }^2}v_2\left( {dz} \right)}  &\le {\alpha '_{3}}{ \vert  y   \vert ^2}+\alpha,	
	\end{align}
	where $\alpha>0$ is a constant, and the inner product is defined in $\mathbb{R}^n$.
	\item For $x \in Q$, a compact subset of ${\mathbb{R}^n}$,  and ${f}\left( {x,\gamma ,\cdot} \right)$ and ${\sigma}\left( {x,\gamma,\cdot} \right)$  can be averaged with respect to the invariant measure $\nu = \left( {{\nu_1},{\nu_2}, \cdots ,{\nu_n}} \right)$ and ${\mu } \left( {dy } \right) $ as follows:
	$$\sum\limits_{\gamma  = 1}^n \int_{\mathbb{R}^n} {{f}\left( {{x},y,\gamma} \right){\mu }} \left( {dy } \right){\nu_\gamma } = {\bar f }\left( {x} \right),$$
	$$\sum\limits_{\gamma  = 1}^n \int_{\mathbb{R}^n} {a_{ij}\left( {{x},y,\gamma } \right){\mu }} \left( {dy } \right) {\nu_\gamma }= \bar a _{ij}\left( {x} \right),$$
	where ${ a_{ij}} = \sum\limits_{k = 1}^n {{\sigma _{ik}}{\sigma _{kj}}} $.
	\item The following averaged equation has a unique weak solution (i.e. uniqueness in the sense of distribution)  on $\left[ {0,T} \right]$,
	\begin{eqnarray}\label{lemma3.2}
	dx\left( t \right) &=& \bar f \left( {x\left( t \right)} \right)dt + \bar \sigma  \left( {x\left( t \right)} \right)dW_1\left( t \right)\cr
	&&+ \int_{ \vert  z  \vert < c} { g \left( {x\left( t \right),z} \right)\tilde N_1\left( {dt,dz} \right)},
	\end{eqnarray}
	where $\sigma=[\sigma _{ij}]_{n \times n}$ and $g=[g _{ij}]_{n \times n}$.
\end{enumerate}	
Throughout this paper, the letter C below with or without subscripts  denotes positive constant whose value may change.		
		
\section{Main results}\label{sec-3}

\begin{thm}\label{thm3.4}
	Let ${x^\varepsilon }\left( t \right)$ be ${\mathbb{R}^n}$-valued and defined on $\left[ { 0,T } \right]$. Suppose the {\rm (A1)} to {\rm (A5)} hold, 
	then,  {${x^\varepsilon }\left(  \cdot  \right)$ converges weakly to  $x\left(  \cdot  \right)$}, where $x\left(  \cdot  \right)$ is the unique weak solution of averaged slow system.
\end{thm}
To prove this theorem, we use the martingale problem formulation. Therefore, 
Define an infinitesimal operator of (\ref{lemma3.2}) as follows
\begin{eqnarray*}
	{\mathcal A}\iota\left( {x} \right) &=& \tilde {\mathcal A}\iota\left( {x} \right) + \mathcal{J }\iota\left( {x} \right) ,
\end{eqnarray*}		
$$\tilde {\mathcal A}\iota\left( {x} \right) = \sum\limits_{i = 1}^{n } {{{\bar f}_i}\left( {x} \right)\frac{{\partial \iota\left( {x} \right)}}{{\partial {x_i}}} + \frac{1}{2}} \sum\limits_{i,j = 1}^{n } {{ {{\bar a}_{ij}}}\left( {x} \right)\frac{{\partial ^2\iota\left( {x} \right)}}{{\partial {x_i}\partial {x_j}}}}, $$
$$\mathcal{J }\iota\left( {x} \right)=\int_{\vert   {z} \vert  < c} \Big[ \iota\left( {x + {g(x, {z})}} \right) - \iota\left( {x} \right) -\sum\limits_{i = 1}^{n } g_i(x, {z}) {\frac{{\partial \iota\left( {x} \right)}}{{\partial {x_i}}} }\Big]v_1(d {z}).$$		
We say that $x\left( t \right)$ solves the martingale problem for the operator $ {\mathcal A}$ if
\begin{eqnarray*}
	{M_\iota} = \iota\left( {x\left( t \right)} \right) - \iota\left( {x\left( 0 \right)} \right) - \int_0^t {{\mathcal A}\iota\left( {x\left( s \right)} \right)}ds
\end{eqnarray*}
is a martingale for each function $\iota$.

The following lemmas are needed.

\begin{lem}\label{lem3.1}\cite{Yin1998Continuous}
	Assume that the generator ${\widetilde {\mathcal{Q}}} $ is weakly irreducible, then for any bounded deterministic function $\beta\left( \cdot \right)$ and each $\gamma \in \mathcal{M}$, when $\varepsilon \to 0$,
	$$\mathbb{E}{\Big \vert  {\int_0^t {\left( {{I_{\left\{ {{r^\varepsilon }\left( u \right) = \gamma } \right\}}} - {\nu _\gamma }} \right){\beta _\gamma }\left( u \right)du} } \Big \vert ^2} \to 0.$$
\end{lem}
\para{Proof:} The proof for this Lemma  refers to \cite{Yin1998Continuous}.
		
\begin{lem}\label{lem3.2}
	Under the conditions {\rm (A1)} to {\rm (A3)}, the slow-fast system $\left\{x^{\varepsilon}(\cdot),y^{\varepsilon}(\cdot)\right\}$ admits a unique solution in ${L^2}\left[{0,T} \right]$ which is a space of square-integrable functions defined on $\left[{0,T} \right]$.
\end{lem}
\para{Proof:} The proof for this Lemma  is similar to \cite{Yin2004Two}, and omitted.

\begin{lem}\label{1lem3.2}
 Under the condition {\rm (A1)} to {\rm (A3)},    
	${y \left( t \right)}$, which is the solution to the following system, 
	$$dy \left( t \right) = {b}\left( {{y}\left( t \right)} \right)dt + V \left( {{y }\left( t \right)} \right)d{ W_2}\left( t \right) + \int_{\vert  {z } \vert < c} {h\left( {{y  }\left( t \right),z} \right){{\tilde N}_{2}}\left( {dt,dz} \right)},$$ 
	where $y(0)=y_0$. Let $\{P_t\}_{t\ge 0}$ be the  transition semigroup of $\{y (t)\}_{t\ge 0}$, i.e., for any bounded continuous function $F: \mathbb{R}^n \to{R}$,
	$P_t F(y):=\mathbb{E} F(y (t)), y\in\mathbb{R}^n, t\ge0. $then $P_t$ has a unique invariant measure $\mu$.
	Moreover, define that ${\bar F }\left( {x} \right)=\int_{\mathbb{R}^n} {{F}\left( {{x},y } \right){\mu }} \left( {dy } \right)$, then
	\begin{eqnarray}\label{lemma3.3}
	\big\vert {P_t F\left( x,{{y}} \right) -\bar F(x) } \big\vert < C{e^{ -  {\lambda t}}}\left( {1 + \vert x \vert + \vert y_0 \vert} \right),
	\end{eqnarray}
	where $\lambda>0$.
\end{lem}

\para{Proof:} The proof is a slight extension of  \cite{Majka2017Coupling} and omitted for the brevity.		
		
\begin{lem}\label{lem3.4}
	If {\rm (A1)} to {\rm (A3)} hold, then  $\left\{ {{x^{\varepsilon} }\left( t \right)} \right\}$ is tight in ${D}({\left[ {0,T } \right]},\mathbb{R}^n) $, the space of $\mathbb{R}^n$-valued functions on the interval $\left[0, T \right]$ that are  {right-continuous} and have left limits.
\end{lem}
\para{Proof:}  Take a slight extension of Lem.2.6 \cite{Yin2004Two} and Lem.2.1 \cite{Givon2007Strong}, it has ,
\begin{eqnarray}\label{lemma3.0}
\mathop {\lim }\limits_{K \to \infty } \mathop {\lim \sup }\limits_{\varepsilon  \to 0} P\Big\{ {\mathop {\sup }\limits_{ 0  \le t \le T} \big\vert  {{x^{\varepsilon }}\left( t \right)} \big\vert  \ge K} \Big\} = 0.
\end{eqnarray}		
For any $\varepsilon  > 0$, there exists a $\tau > 0$, with  aid of the H\"older inequality and the Burkholder-Davis-Gundy inequality in \cite{Protter2004Stochastic},
\begin{eqnarray*}\label{lemma3.5}
	\mathbb{E}\big[ {  {{\vert   {{x^{\varepsilon }}\left( {t + \tau } \right) - {x^{\varepsilon }}\left( t \right)} \vert }^2}} \big]
	&\le& C  \mathbb{E}{\Big[ {\int_t^{t + \tau } {f\left( {{x^{\varepsilon }}\left( s \right),{y ^\varepsilon }\left( s \right),{r^\varepsilon }\left( s \right)} \right)} ds} \Big]^2}\cr
	&&+ C  \mathbb{E}{\Big[ { \int_t^{t + \tau } {\sigma \left( {{x^{\varepsilon }}\left( s \right),{y ^\varepsilon }\left( s\right),{r^\varepsilon }\left( s \right)} \right)} dW_1\left( s \right)} \Big]^2}\cr
	&&+ C  \mathbb{E}\Big[  \int_t^{t + \tau } \int_{\vert  z \vert  < c} g\left( {{x^{\varepsilon }}\left( s\right),z} \right)\tilde N_1\left( {dz,ds} \right)   \Big]^2\cr
	&\le&   {C'_\gamma} \mathbb{E}\Big[ \int_t^{t + \tau } \big( 1 + {{\big\vert  {{x^{\varepsilon }}\left( s \right)} \big\vert }^2}+ {{\big\vert  {{y ^\varepsilon }\left( s \right)} \big\vert }^2} \big) ds \Big].
\end{eqnarray*}
Take a slight extension of Lem.2.1 \cite{Givon2007Strong}, we obtain the boundness of $y^\varepsilon(\cdot)$ as follows
\begin{eqnarray*}\label{lemma3.6}
	\mathbb{E}\big[ {{{\vert  {y ^\varepsilon \left( t \right)} \vert }^2}} \big] \le {C },
\end{eqnarray*}		
where ${C }$ is a constant which is dependent on the initial data of $y ^\varepsilon\left( t \right)$.
Together with the boundness of $ {{x^\varepsilon }\left( t \right)} $ in probability shown in \lemref{lem3.2}, it yields that
\begin{eqnarray*}\label{lemma3.7}
	\mathbb{E}\big[ { {{\vert {{x^{\varepsilon }}\left( {t + \tau } \right) - {x^{\varepsilon }}\left( t \right)} \vert }^2}} \big]
	&\le& {C_T} \cdot {C_K} \cdot \tau  =  {C\cdot\tau }.
\end{eqnarray*}		
	By virtue of the criterion  [\cite{Kushner1984Approximation},Theorem 3, p.47], $\{ {{x^\varepsilon }\left(  \cdot  \right)} \}$ is tight in ${D} {(\left[ {0,T } \right],\mathbb{R}^n)} $.	
		
Since $x^\varepsilon(\cdot)$ is tight, by Prohorov's theorem, it is sequentially compact. Thus we can extract a weakly convergent subsequence. Do so and still denote the convergent subsequence by $\varepsilon$ and denote the limit as $x(\cdot)$. To characterize the limit process $x(\cdot)$, the martingale method  is shown in the following lemma:
\begin{lem}\label{thm3.5}
	Let ${x^\varepsilon }\left( t \right)$ be ${\mathbb{R}^n}$-valued and defined on $\left[ { 0 ,T } \right)$. If $\left\{ {{x^\varepsilon }\left(  \cdot  \right)} \right\}$ be tight on ${D}({\left[ {0,T } \right]},\mathbb{R}^n)$,  and for each $\iota \left(  \cdot  \right)\in C_0^4\left( {{\mathbb{R}^n},\mathbb{R}} \right)$, set of $C^4$ functions with
	compact support, and each $T < \infty $, then there exists a $\iota^\varepsilon\left(  \cdot  \right) \in D( {{{\hat {\mathcal A}}^\varepsilon }} )$,  the domain of ${\hat {\mathcal A}^\varepsilon }$, such that
	\begin{eqnarray}\label{lemma3.9}
	p \text{-}{ \mathop {\lim }\limits_{\varepsilon  \to 0} \big[ {{\iota^\varepsilon }\left(  \cdot \right) - \iota\left( {{x^\varepsilon }\left(  \cdot  \right)} \right)} \big]} = 0,
	\end{eqnarray}
	and
	\begin{eqnarray}\label{lemma3.10}
	\mathop {\lim }\limits_{\varepsilon  \to 0}\mathbb{E}\big\vert \int_{ {s}}^{T} \big[{{\hat {\mathcal A}}^\varepsilon }{\iota^{\varepsilon }}\left( t \right)-{\mathcal A}\iota\left( {{x^{\varepsilon }}\left( t \right)} \right)\big]dt\big\vert =0.
	\end{eqnarray}
	for each $t<T$.
	Then, ${x^\varepsilon }\left(  \cdot  \right)$ weakly converges to  $x\left(  \cdot  \right)$, where $x\left(  \cdot  \right)$ is the unique weak solution of (\ref{lemma3.2}).
\end{lem}
\para{Proof:} This is a slight extension of Theorem 2 in \cite{Kushner1984Approximation}.		
		
\para{Proof of  \thmref{thm3.4}:} 	
		 {Step 1: Choose appropriate perturbed test functions and verify the Condition (\ref{lemma3.9}).} 
		
		According to the definition of $p-lim$, to prove (\ref{lemma3.9}) for $x^\varepsilon(t)$, for each $\iota \left(  \cdot  \right)\in C_0^4\left( {{\mathbb{R}^n},\mathbb{R}} \right)$,  and each $T < \infty $, we need to look for a function $\iota^\varepsilon\left(  \cdot  \right) \in D( {{{\hat {\mathcal A}}^\varepsilon }} )$,  and verify the corresponding conditions,
		\begin{eqnarray}\label{a1}
			\begin{split}
				\begin{cases}
					\mathop {\sup }\limits_{t,\varepsilon } \mathbb{E} \vert  {{\iota^\varepsilon }\left(  t \right) - \iota\left( {{x^\varepsilon }\left(  t  \right)} \right)} \vert  < \infty ,\cr
					\mathop {\lim }\limits_{\varepsilon  \to 0} \mathbb{E} \vert  {{\iota^\varepsilon }\left(  t \right) - \iota\left( {{x^\varepsilon }\left(  t  \right)} \right)}  \vert  = 0 ,
				\end{cases} 
			\end{split}
		\end{eqnarray}
		for each $t$. 

Next, we use $x$ to represent ${x^{\varepsilon }}\left( t \right)$. 
For any $\iota \in C_0^4\left( {{\mathbb{R}^n},\mathbb{R}} \right)$, we  define the following perturbed test functions,
\begin{eqnarray}\label{lemma3.12}
{\iota^{\varepsilon }}\left( t \right) = \iota( {{x}} ) + \iota_1^{\varepsilon }\left( t \right) + \frac{1}{2}\iota_2^{\varepsilon }\left( t \right) ,
\end{eqnarray}		
where
$$\iota_1^{\varepsilon }\left( t \right) = \sum\limits_{i = 1}^{n } {\int_t^T {{\iota_{{x_i}}}\left( x \right){\mathbb{E}_t^\varepsilon}\big[ {\bar f^\nu_i\left( {x,{y ^\varepsilon }\left( {u } \right)} \right) - \bar f_i\left( {x} \right)} \big]du} },$$
$$\iota_2^{\varepsilon }\left( t \right) = \sum\limits_{i,j = 1}^{n } {\int_t^T {{\iota_{{x_i}{x_j}}}\left( x \right){\mathbb{E}_t^\varepsilon}\big[ {\bar a^\nu_{ij}\left( {x,{y ^\varepsilon }\left( {u} \right)} \right) - \bar a_{ij}\left( {x} \right)} \big]du} },$$		
and 
$$\bar f^\nu_i\left( {x,{y  }} \right)=\sum\limits_{\gamma  = 1}^n  {{f_i}\left( {{x},y,\gamma } \right)} {\nu_\gamma },\bar a_{ij}^\nu\left( {x,{y  }} \right)=\sum\limits_{\gamma  = 1}^n  {{a_{ij}}\left( {{x},y,\gamma } \right)} {\nu_\gamma },$$
where $\iota_{x_i}(x)=\frac{\partial {\iota(x)}}{\partial{x_i}}$ for $i\in\{1,2,...,n\}$ and $\iota_{x_i x_j}(x)=\frac{\partial^2 {\iota(x)}}{\partial{x_i}\partial{x_j}}$ for $i,j\in\{1,2,...,n\}$, $y^\varepsilon(u)=y(u / \varepsilon)$.
We also note that for any bounded function $F\left( {u \ge t} \right)$, owing to the independence between the Markovian switching and fast component,	
\begin{eqnarray}\label{lemma3.11}
{\mathbb{E}_t^\varepsilon}F\left( {{y ^\varepsilon }\left( {u,{y ^\varepsilon }\left( t \right)} \right)} \right) =\int_{{\mathbb{R}^n}} F\left( \xi \right){P}\left( {{y ^\varepsilon }\left( t \right),t,d\xi,u} \right),
\end{eqnarray}
where ${P}\left( {y,t,A,u} \right) = P\left( {{y ^\varepsilon }\left( {u} \right) \in A \vert  {y ^\varepsilon }\left( {t} \right) = y} \right)$.	
	
Firstly we estimate $ \iota_1^{\varepsilon }\left( t \right)$,  $\iota_2^{\varepsilon }\left( t \right)$ and $\iota_3^{\varepsilon }\left( t \right)$ respectively.
With the \lemref{1lem3.2}, we have
\begin{eqnarray}\label{lemma3.14}
\mathop {\lim }\limits_{\varepsilon  \to 0} \mathbb{E}  \big\vert  {\iota_1^{\varepsilon }\left( t \right)} \big\vert  &\le&\mathop {\lim }\limits_{\varepsilon  \to 0} \mathbb{E}\Big\{ \varepsilon \sum\limits_{i = 1}^{n } {  \Big\vert  {\int_{{t \mathord{\left/{\vphantom {t \varepsilon }} \right.\kern-\nulldelimiterspace} \varepsilon }}^{{T \mathord{\left/{\vphantom {T \varepsilon }} \right.\kern-\nulldelimiterspace} \varepsilon }} {{\iota_{{x_i}}}\left( {x} \right) {\mathbb{E}_{t/ \varepsilon}^\varepsilon}\big[ {\bar f^\nu_i\left( {{x},{y  }\left( {u} \right)} \right) - {\bar {f_i}}\left( {{x}} \right)} \big]du} } \Big\vert } \Big\}\cr
&\le& \mathop {\lim }\limits_{\varepsilon  \to 0} \mathbb{E}\Big\{\varepsilon {C_1}\sum\limits_{i = 1}^{n } {  \Big\vert  {\int_{{t \mathord{\left/{\vphantom {t \varepsilon }} \right.\kern-\nulldelimiterspace} \varepsilon }}^{{T \mathord{\left/{\vphantom {T \varepsilon }} \right.\kern-\nulldelimiterspace} \varepsilon }} {\left\vert  {{\iota_{{x_i}}}\left( {x} \right)} \right\vert {e^{ - \frac{{\lambda  (u-t)}}{2}}}du} } \Big\vert } \Big\}\cr
&=& \mathop {\lim }\limits_{\varepsilon  \to 0}\sum\limits_{i = 1}^{n }  \mathbb{E}\Big\{\frac{{{\varepsilon} }}{\lambda } \cdot C \cdot \big\vert  {{\iota_{{x_i}}}\left( x \right)} \big\vert  \cdot \big( {1 - {e^{ - \frac{{\lambda \left(T-t\right)}}{{ {2}\varepsilon }}}}} \big)\Big\}\cr
&=& 0.
\end{eqnarray}

Similarly, we obtain that
\begin{eqnarray}\label{lemma3.16}
\mathop {\lim }\limits_{\varepsilon  \to 0} \mathbb{E}  \big\vert  {\iota_2^{\varepsilon }\left( t \right)} \big\vert 
=0.
\end{eqnarray}

From (\ref{lemma3.14}) to (\ref{lemma3.16}), we prove that (\ref{lemma3.12}) satisfies the condition (\ref{a1}). According to the definition of $p-lim$, that is, (\ref{lemma3.9}) holds.	
	
 {Step 2: Then we prove that $\iota (x^\varepsilon (t))$ also satisfies the condition (\ref{lemma3.10}).} 

To proceed, ${\hat {\mathcal A}^\varepsilon }\iota_1^{\varepsilon }\left(  \cdot  \right)$ is estimated as follows:
\begin{eqnarray}\label{lemma3.18}
\hat{\mathcal A}^\varepsilon \iota_1^{\varepsilon }\left( t \right) &=& p\text{-}{\mathop {\lim }\limits_{\delta\to 0}{\frac{{\mathbb{E}_t^\varepsilon \iota_1^{\varepsilon }\left({t+\delta}\right)-\iota_1^{\varepsilon }\left(t\right)}}{\delta}}} \cr
&=:& \Xi_1+\Xi_2,
\end{eqnarray}
where
\begin{eqnarray*}
	{\Xi_1} &=&  - \sum\limits_{i = 1}^{n }{\iota_{{x_i}}}\left( x \right)\big[ {\bar f^\nu_i\left( {x,{y ^\varepsilon }\left( {t} \right)} \right) - \bar f_i\left( {x} \right)} \big]\cr
	&&+ \sum\limits_{i,j = 1}^{n } {\varepsilon \int_{{t \mathord{\left/{\vphantom {t \varepsilon }} \right.\kern-\nulldelimiterspace} \varepsilon }}^{{T \mathord{\left/{\vphantom {T \varepsilon }} \right.\kern-\nulldelimiterspace} \varepsilon }} {{ {\iota_{{x_i}{x_j}}}}\left( x \right){ {\mathbb{E}_{t/ \varepsilon}^\varepsilon}}\big[ {\bar f^\nu_i\left( {x,{y ^\varepsilon }\left( {u} \right)} \right) - \bar f_i\left( {x} \right)} \big]du} } f_j\left( {x,{y ^\varepsilon  }\left( t \right),r^\varepsilon(t)} \right)\cr
	&&+ \sum\limits_{i,j = 1}^{n } {\varepsilon \int_{{t \mathord{\left/{\vphantom {t \varepsilon }} \right.\kern-\nulldelimiterspace} \varepsilon }}^{{T \mathord{\left/{\vphantom {T \varepsilon }} \right.\kern-\nulldelimiterspace} \varepsilon }} {{\iota_{{x_i}}}\left( x \right){ {\mathbb{E}_{t/ \varepsilon}^\varepsilon}}{{\big[ {\bar f^\nu_i\left( {x,{y  }\left( {u} \right)} \right) - \bar f_i\left( {x} \right)} \big]}_{{x_j}}}du} } f_j\left( {x,{y ^\varepsilon  }\left( t \right),r^\varepsilon(t)} \right)\cr
	&&+ \frac{1}{2}\sum\limits_{i,j,k = 1}^{n } {\varepsilon \int_{{t \mathord{\left/{\vphantom {t \varepsilon }} \right.\kern-\nulldelimiterspace} \varepsilon }}^{{T \mathord{\left/{\vphantom {T \varepsilon }} \right.\kern-\nulldelimiterspace} \varepsilon }} {{\iota_{{x_i}{x_j}{x_k}}}\left( x \right){ {\mathbb{E}_{t/ \varepsilon}^\varepsilon}}\big[ {\bar f^\nu_i\left( {x,{y  }\left( {u} \right)} \right) - \bar f_i\left( {x} \right)} \big]du} } a_{jk}\left( {x,{y ^\varepsilon  }\left( t \right),r^\varepsilon(t)} \right)\cr
	&&+\sum\limits_{i,j,k = 1}^{n } {\varepsilon \int_{{t \mathord{\left/{\vphantom {t \varepsilon }} \right.\kern-\nulldelimiterspace} \varepsilon }}^{{T \mathord{\left/{\vphantom {T \varepsilon }} \right.\kern-\nulldelimiterspace} \varepsilon }} {{\iota_{{x_i}{x_j}}}\left( x \right){ {\mathbb{E}_{t/ \varepsilon}^\varepsilon}}{{\big[ {\bar f^\nu_i\left( {x,{y  }\left( {u} \right)} \right) - \bar f_i\left( {x} \right)} \big]}_{{x_k}}}du} } a_{jk}\left( {x,{y ^\varepsilon  }\left( t \right),r^\varepsilon(t)} \right)\cr
	&&+ \frac{1}{2}\sum\limits_{i,j,k = 1}^{n } {\varepsilon \int_{{t \mathord{\left/{\vphantom {t \varepsilon }} \right.\kern-\nulldelimiterspace} \varepsilon }}^{{T \mathord{\left/{\vphantom {T \varepsilon }} \right.\kern-\nulldelimiterspace} \varepsilon }} {{\iota_{{x_i}}}\left( x \right){ {\mathbb{E}_{t/ \varepsilon}^\varepsilon}}{{\big[ {\bar f^\nu_i\left( {x,{y  }\left( {u} \right)} \right) - \bar f_i\left( {x} \right)} \big]}_{{x_j}{x_k}}}du} } a_{jk}\left( {x,{y ^\varepsilon  }\left( t \right),r^\varepsilon(t)} \right)\cr
	&=:& \sum\limits_{m = 1}^{6} {\Xi_1^{\left( m \right)}},
\end{eqnarray*}
 {where ${\Xi_1^{\left( 1 \right)}}$ is from the differential on $t$ in this integral, and ${\Xi_1^{\left( m \right)}}$ with $m\in \{2,3\}$ are from the differentials on $x^\varepsilon(t)$, and ${\Xi_1^{\left( m \right)}}$ with $m\in \{4,5,6\}$ are from the second-order differentials on $x^\varepsilon(t)$. And ${\Xi_2}$ is from the difference on jump of $x^\varepsilon(t)$.} Take the Taylor expansion, then 
\begin{eqnarray*}
	{\Xi_2} 
	&=& \frac{1}{2}\int_{\vert { {z}} \vert  < c} \Big\{ \sum\limits_{i,j,k = 1}^{n } \varepsilon \int_{{t \mathord{\left/{\vphantom {t \varepsilon }} \right.\kern-\nulldelimiterspace} \varepsilon }}^{{T \mathord{\left/{\vphantom {T \varepsilon }} \right.\kern-\nulldelimiterspace} \varepsilon }} {\iota_{{x_i}{x_j}{x_k}}}\left( x+ {\theta g(x, {z})} \right){ {\mathbb{E}_{t/ \varepsilon}^\varepsilon}}\big[ \bar f^\nu_i\left( {x+ {\theta g(x, {z})},{y  }\left( {u} \right)} \right)\big.\Big.\cr
	&&\qquad\Big. \big.- \bar f_i\left( {x+ {\theta g(x, {z})}} \right) \big]du   \Big. \cr
	&&+ 2\sum\limits_{i,j,k = 1}^{n } \varepsilon \int_{{t \mathord{\left/{\vphantom {t \varepsilon }} \right.\kern-\nulldelimiterspace} \varepsilon }}^{{T \mathord{\left/{\vphantom {T \varepsilon }} \right.\kern-\nulldelimiterspace} \varepsilon }} {\iota_{{x_i}{x_j}}}\left( x+ {\theta g(x, {z})} \right){ {\mathbb{E}_{t/ \varepsilon}^\varepsilon}}\big[ \bar f^\nu_i\left( {x+ {\theta g(x, {z})},{y  }\left( {u} \right)} \right) \big.\Big.\cr
	&&\qquad\Big. \big.- \bar f_i\left( {x+ {\theta g(x, {z})}} \right) \big]_{{x_k}}du \Big. \cr
	&&+ \sum\limits_{i,j,k = 1}^{n } \varepsilon \int_{{t \mathord{\left/{\vphantom {t \varepsilon }} \right.\kern-\nulldelimiterspace} \varepsilon }}^{{T \mathord{\left/{\vphantom {T \varepsilon }} \right.\kern-\nulldelimiterspace} \varepsilon }} {\iota_{{x_i}}}\left( x+ {\theta g(x, {z})} \right){ {\mathbb{E}_{t/ \varepsilon}^\varepsilon}}\big[ \bar f^\nu_i\left( {x+ {\theta g(x, {z})},{y  }\left( {u} \right)} \right)\big.\Big.\cr
	&&\qquad\Big. \big.- \bar f_i\left( {x+ {\theta g(x, {z})}} \right) \big]_{{x_j}{x_k}}du \Big\}G_{jk}v\left( {d {z}} \right) \cr
	&=:& \sum\limits_{m = 1}^3 {\Xi_2^{\left( m \right)}},
\end{eqnarray*}
where $g ={g(x, {z})}$ and  ${ G_{jk}} = \sum\limits_{h = 1}^n {{g _{jh}}{g_{hk}}} $ and  {$0<\theta<1$}.	
	
Then we compute $\Xi_1^{\left( m \right)}, m=1,\ldots,7$ respectively.	
Taking the similar  argument to (\ref{lemma3.14}), one has	
\begin{eqnarray}\label{lemma3.20}
 {\Xi_1^{\left( 2 \right)}} &=&   \varepsilon \sum\limits_{i,j = 1}^{n }  \  \int_{{t \mathord{\left/{\vphantom {t \varepsilon }} \right.\kern-\nulldelimiterspace} \varepsilon }}^{{T \mathord{\left/{\vphantom {T \varepsilon }} \right.\kern-\nulldelimiterspace} \varepsilon }} {\iota_{{x_i}{x_j}}}\left( x \right){ {\mathbb{E}_{t/ \varepsilon}^\varepsilon}}\big[ {\bar f^\nu_i\left( {x,{y  }\left( {u} \right)} \right) - \bar f_i\left( {x} \right)} \big]duf_j\left( {x,{y ^\varepsilon  }\left( t \right),r^\varepsilon(t)} \right)   \cr
&\le&  \mathbb{E}\varepsilon {C_1}\sum\limits_{i,j = 1}^{n } {  {\int_{{t \mathord{\left/{\vphantom {t \varepsilon }} \right.\kern-\nulldelimiterspace} \varepsilon }}^{{T \mathord{\left/{\vphantom {T \varepsilon }} \right.\kern-\nulldelimiterspace} \varepsilon }} { \big\vert {{\iota_{{x_i}{x_j}}}\left( {x} \right)} {e^{ - \frac{{\lambda  (u-t)}}{ {2}}}}du\big\vert f_j\left( {{x},{r^{\varepsilon }( t)},{y ^\varepsilon }\left( t\right)} \right)} } } \cr
&=&  \sum\limits_{i,j = 1}^{n }\mathbb{E}\Big\{\frac{{{\varepsilon} }}{\lambda } \cdot C \cdot \big\vert  {{\iota_{{x_i}{x_j}}}\left( {x} \right)} \big\vert  \cdot \big( {1 - {e^{ - \frac{{\lambda \left(T-t\right)}}{{ {2}\varepsilon }}}}} \big)\Big\}f_j\left( {{x},{y ^\varepsilon  }\left( t \right),r^\varepsilon(t)} \right)\cr
&=& O(\varepsilon).
\end{eqnarray}
Similarly,
\begin{eqnarray}\label{lemma3.23}
  \mathbb{E}\big\vert {\Xi_1^{\left( 4 \right)}}\big\vert  =O(\varepsilon).
\end{eqnarray}
Then we estimate ${\Xi_1^{\left( 3 \right)}}$ as follows. Owing to the fact that the fast component is independent to the slow component, by the \lemref{lem3.1},   
\begin{eqnarray}\label{lemma3.24}
  {\Xi_1^{\left( 3 \right)}} &=&   \varepsilon \sum\limits_{i,j = 1}^{n }  \  \int_{{t \mathord{\left/{\vphantom {t \varepsilon }} \right.\kern-\nulldelimiterspace} \varepsilon }}^{{T \mathord{\left/{\vphantom {T \varepsilon }} \right.\kern-\nulldelimiterspace} \varepsilon }} { {\mathbb{E}_{t/ \varepsilon}^\varepsilon}}\big[ {\bar f^\nu_i\left( {x,{y  }\left( {u} \right)} \right) - \bar f_i\left( {x} \right)} \big]_{x_j}duf_j\left( {x,{y ^\varepsilon  }\left( t \right),r^\varepsilon(t)} \right)   \cr
&=&  \Big\{\varepsilon \sum\limits_{i,j = 1}^{n }   \int_{{t \mathord{\left/{\vphantom {t \varepsilon }} \right.\kern-\nulldelimiterspace} \varepsilon }}^{{T \mathord{\left/{\vphantom {T \varepsilon }} \right.\kern-\nulldelimiterspace} \varepsilon }} {\iota_{{x_i}}}\left({x} \right){ {\mathbb{E}_{t/ \varepsilon}^\varepsilon}}\big[ {\bar f^\nu_{i,{x_j}}\left( {{x},{y  }\left( {u} \right)} \right) - { {\bar f_{i,{x_j}}}}\left( {{x}} \right)} \big] duf_j\left( {{x},{y ^\varepsilon  }\left( t \right),r^\varepsilon(t)} \right)\cr
&\le&  {\varepsilon} {C_1}\sum\limits_{i,j = 1}^{n }\mathbb{E}  {\Big\{{  {\int_{{t \mathord{\left/{\vphantom {t \varepsilon }} \right.\kern-\nulldelimiterspace} \varepsilon }}^{{T \mathord{\left/{\vphantom {T \varepsilon }} \right.\kern-\nulldelimiterspace} \varepsilon }} { {{\iota_{{x_i}{x_j}}}\left( {x} \right)} {e^{ - \frac{{\lambda  (u-t)}}{ {2}}}}du} } } \Big\}}f_j\left( {x,{y ^\varepsilon  }\left( t \right),r^\varepsilon(t)} \right)\cr
&=&O(\varepsilon).
\end{eqnarray}
Therefore we get
\begin{eqnarray*}\label{lemma3.27}
	 \mathbb{E}\big\vert {\Xi_1^{\left( 3 \right)}}\big\vert  =O(\varepsilon).
\end{eqnarray*}
And taking the similar methods, we arrive at
\begin{eqnarray*}\label{lemma3.28}
  \mathbb{E}\big\vert {\Xi_1^{\left( 5 \right)}}\big\vert  =O(\varepsilon),\quad \mathbb{E}\big\vert {\Xi_1^{\left( 6 \right)}}\big\vert =O(\varepsilon).
\end{eqnarray*}	
Then we estimate ${\Xi_2}$, firstly  ${\Xi_2^{\left( 1 \right)}}$ is computed by the similar  argument to the ${\Xi_1^{\left( 2 \right)}}$,	
\begin{eqnarray}\label{lemma3.31}
  \mathbb{E}\big\vert {\Xi_2^{\left( 1 \right)}}\big\vert  &=&  \frac{1}{2}\sum\limits_{i,j,k = 1}^{n }\mathbb{E}\Big\vert  \int_{\vert  { {z}} \vert  < c}  { {\varepsilon \int_{{t \mathord{\left/{\vphantom {t \varepsilon }} \right.\kern-\nulldelimiterspace} \varepsilon }}^{{T \mathord{\left/{\vphantom {T \varepsilon }} \right.\kern-\nulldelimiterspace} \varepsilon }} {{\iota_{{x_i}{x_j}{x_k}}}\left( x \right){ {\mathbb{E}_{t/ \varepsilon}^\varepsilon}}\big[ {\bar f^\nu_i\left( {x,{y  }\left( {u} \right)} \right) - \bar f_i\left( {x} \right)} \big]du} } }\Big.\cr
&&\quad\Big. {G_{jk}}v\left( {d {z}} \right) \Big\vert \cr
&\le&  \frac{1}{2} \varepsilon\sum\limits_{i ,j,k= 1}^{n }\mathbb{E}\Big\vert \int_{\vert  { {z}} \vert  < c} \Big\{ {C_1}   \int_{{t \mathord{\left/{\vphantom {t \varepsilon }} \right.\kern-\nulldelimiterspace} \varepsilon }}^{{T \mathord{\left/{\vphantom {T \varepsilon }} \right.\kern-\nulldelimiterspace} \varepsilon }} \vert  {{\iota_{{x_i}{x_j}{x_k}}}\left( {x} \right)} \vert {e^{ - \frac{{\lambda  (u-t)}}{ {2}}}}du   \Big\} {G_{jk}}v\left( {d {z}} \right) \Big\vert \cr
&=& \frac{1}{2} C{\varepsilon} \sum\limits_{i,j,k = 1}^{n }\int_{\vert  { {z}} \vert  <c} {\big\vert  {{\iota_{{x_i}{x_j}{x_k}}}\left( {x} \right)} \big\vert \big( {1 - {e^{ - \frac{{\lambda \left(T-t\right)}}{{ {2}\varepsilon }}}}} \big){ {G_{jk}}}v\left( {d {z}} \right)} \cr
&=& O(\varepsilon).
\end{eqnarray}
Similarly, taking the estimation of ${\Xi_2^{\left( 2 \right)}}$ and ${\Xi_2^{\left( 3 \right)}}$, it yields that
\begin{eqnarray}\label{lemma3.32}
  \mathbb{E}\big\vert {\Xi_2^{\left( 2 \right)}}\big\vert =O(\varepsilon),\quad   \mathbb{E}\big\vert {\Xi_2^{\left( 3 \right)}}\big\vert =O(\varepsilon) .
\end{eqnarray}
From the (\ref{lemma3.31}) to (\ref{lemma3.32}),
\begin{eqnarray}\label{lemma3.33}
 \mathbb{E}\big\vert {\Xi_2}\big\vert =O(\varepsilon) .
\end{eqnarray}
Together with (\ref{lemma3.18}) and (\ref{lemma3.33}), it has
\begin{eqnarray}\label{lemma3.35}
\hat{\mathcal A}^\varepsilon \iota_1^{\varepsilon }\left( t \right) =  - \sum\limits_{i = 1}^{n }{\iota_{{x_i}}}\left( x \right) {\bar f^\nu_i\left( {x,{y ^\varepsilon }\left( {t} \right)} \right) + \sum\limits_{i = 1}^{n }\iota_{{x_i}} \bar f_i\left( {x} \right)}+O(\varepsilon).
\end{eqnarray}
Then estimate $\hat{\mathcal A}^\varepsilon \iota_2^{\varepsilon }\left( t \right) $ as following:
\begin{eqnarray*}\label{lemma3.36}
	\hat{\mathcal A}^\varepsilon \iota_2^{\varepsilon }\left( t \right) &=& p\text{-}{\mathop {\lim }\limits_{\delta\to 0}{\frac{{\mathbb{E}_t^\varepsilon \iota_2^{\varepsilon }\left({t+\delta}\right)-\iota_2^{\varepsilon }\left(t\right)}}{\delta}}} \cr
	&=:& \Xi_3+\Xi_4,
\end{eqnarray*}
where
\begin{eqnarray*}
	{\Xi_3} &=&  - \sum\limits_{i,j = 1}^{n }{\iota_{{x_i}{x_j}}}\left( x\right)\big[ {\bar a^\nu_{ij}\left( {x,{y ^\varepsilon }\left( {t} \right)} \right) - \bar a_{ij}\left( {x} \right)} \big]\cr
	&&+ \sum\limits_{i,j,k = 1}^{n } {\varepsilon \int_{{t \mathord{\left/{\vphantom {t \varepsilon }} \right.\kern-\nulldelimiterspace} \varepsilon }}^{{T \mathord{\left/{\vphantom {T \varepsilon }} \right.\kern-\nulldelimiterspace} \varepsilon }} {{\iota_{{x_i}{x_j}{x_k}}}\left( x \right){ {\mathbb{E}_{t/ \varepsilon}^\varepsilon}}\big[ {\bar a^\nu_{ij}\left( {x,{y  }\left( {u} \right)} \right) - \bar a_{ij}\left( {x} \right)} \big]du} } f_k\left( {x,{y ^\varepsilon  }\left( t \right),r^\varepsilon(t)} \right)\cr
	&&+ \sum\limits_{i,j,k = 1}^{n } {\varepsilon \int_{{t \mathord{\left/{\vphantom {t \varepsilon }} \right.\kern-\nulldelimiterspace} \varepsilon }}^{{T \mathord{\left/{\vphantom {T \varepsilon }} \right.\kern-\nulldelimiterspace} \varepsilon }} {{\iota_{{x_i}{x_j}}}\left( x \right){ {\mathbb{E}_{t/ \varepsilon}^\varepsilon}}{{\big[ {\bar a^\nu_{ij}\left( {x,{y  }\left( {u} \right)} \right) - \bar a_{ij}\left( {x} \right)} \big]}_{{x_k}}}du} } f_k\left( {x,{y ^\varepsilon  }\left( t \right),r^\varepsilon(t)} \right)\cr
	&&+ \frac{1}{2}\sum\limits_{i,j,k,l = 1}^{n } {\varepsilon \int_{{t \mathord{\left/{\vphantom {t \varepsilon }} \right.\kern-\nulldelimiterspace} \varepsilon }}^{{T \mathord{\left/{\vphantom {T \varepsilon }} \right.\kern-\nulldelimiterspace} \varepsilon }} {{\iota_{{x_i}{x_j}{x_k}{x_l}}}\left( x \right){ {\mathbb{E}_{t/ \varepsilon}^\varepsilon}}\big[ {\bar a^\nu_{ij}\left( {x,{y  }\left( {u} \right)} \right) - \bar a_{ij}\left( {x} \right)} \big]du} } a_{kl}\left( {x,{y ^\varepsilon  }\left( t \right),r^\varepsilon(t)} \right)\cr
	&&+\sum\limits_{i,j,k,l = 1}^{n } {\varepsilon \int_{{t \mathord{\left/{\vphantom {t \varepsilon }} \right.\kern-\nulldelimiterspace} \varepsilon }}^{{T \mathord{\left/{\vphantom {T \varepsilon }} \right.\kern-\nulldelimiterspace} \varepsilon }} {{\iota_{{x_i}{x_j}{x_k}}}\left( x \right){ {\mathbb{E}_{t/ \varepsilon}^\varepsilon}}{{\big[ {\bar a^\nu_{ij}\left( {x,{y  }\left( {u} \right)} \right) - \bar a_{ij}\left( {x} \right)} \big]}_{{x_l}}}du} } a_{kl}\left( {x,{y ^\varepsilon  }\left( t \right),r^\varepsilon(t)} \right)\cr
	&&+ \frac{1}{2}\sum\limits_{i,j,k,l = 1}^{n } {\varepsilon \int_{{t \mathord{\left/{\vphantom {t \varepsilon }} \right.\kern-\nulldelimiterspace} \varepsilon }}^{{T \mathord{\left/{\vphantom {T \varepsilon }} \right.\kern-\nulldelimiterspace} \varepsilon }} {{\iota_{{x_i}{x_j}}}\left( x \right){ {\mathbb{E}_{t/ \varepsilon}^\varepsilon}}{{\big[ {\bar a^\nu_{ij}\left( {x,{y  }\left( {u} \right)} \right) - \bar a_{ij}\left( {x} \right)} \big]}_{{x_k}{x_l}}}du} } a_{kl}\left( {x,{y ^\varepsilon  }\left( t \right),r^\varepsilon(t)} \right)\cr
	&=:& \sum\limits_{m = 1}^{6} {\Xi_3^{\left(m \right)}},
\end{eqnarray*}
and take the  {Taylor} expansion, then 
\begin{eqnarray*}
	{\Xi_4}
	&=& \frac{1}{2}\int_{\vert  { {z}} \vert  < c} \Big\{ \sum\limits_{i,j,k,l = 1}^{n } \varepsilon \int_{{t \mathord{\left/{\vphantom {t \varepsilon }} \right.\kern-\nulldelimiterspace} \varepsilon }}^{{T \mathord{\left/{\vphantom {T \varepsilon }} \right.\kern-\nulldelimiterspace} \varepsilon }} {\iota_{{x_i}{x_j}{x_k}{x_l}}}\left( x+ {\theta g(x, {z})} \right){ {\mathbb{E}_{t/ \varepsilon}^\varepsilon}}\big[ \bar a^\nu_{ij}\left( {x+ {\theta g(x, {z})},{y  }\left( {u} \right)} \right)\big.\Big.\cr
	&&\qquad\Big.\big.- \bar a_{ij}\left( {x+ {\theta g(x, {z})}} \right) \big]du   \Big. \cr
	&&+ 2\sum\limits_{i,j,k,l = 1}^{n } \varepsilon \int_{{t \mathord{\left/{\vphantom {t \varepsilon }} \right.\kern-\nulldelimiterspace} \varepsilon }}^{{T \mathord{\left/{\vphantom {T \varepsilon }} \right.\kern-\nulldelimiterspace} \varepsilon }} {\iota_{{x_i}{x_j}{x_k}}}\left( x+ {\theta g(x, {z})} \right){ {\mathbb{E}_{t/ \varepsilon}^\varepsilon}}\big[ \bar a^\nu_{ij}\left( {x+ {\theta g(x, {z})},{y  }\left( {u} \right)} \right)\big.\Big.\cr
	&&\qquad\Big.\big.- \bar a_{ij}\left( {x+ {\theta g(x, {z})}} \right) \big]_{{x_l}}du  \cr
	&&+ \sum\limits_{i,j,k,l = 1}^{n } \varepsilon \int_{{t \mathord{\left/{\vphantom {t \varepsilon }} \right.\kern-\nulldelimiterspace} \varepsilon }}^{{T \mathord{\left/{\vphantom {T \varepsilon }} \right.\kern-\nulldelimiterspace} \varepsilon }} {\iota_{{x_i}{x_j}}}\left( x+ {\theta g(x, {z})} \right){ {\mathbb{E}_{t/ \varepsilon}^\varepsilon}}\big[ \bar a^\nu_{ij}\left( {x+ {\theta g(x, {z})},{y  }\left( {u} \right)} \right)\big.\Big.\cr
	&&\qquad\Big.\big. - \bar a_{ij}\left( {x+ {\theta g(x, {z})}} \right) \big]_{{x_k}{x_l}}du \Big\} {G}_{kl}v\left( {d {z}} \right)\cr
	&=:& \sum\limits_{m= 1}^3 {\Xi_4^{\left( m \right)}},
\end{eqnarray*}
where  {$0<\theta<1$}.
Therefore from  (\ref{lemma3.18}) to (\ref{lemma3.35}), one can obtain the similar conclusion for $\hat{\mathcal A}^\varepsilon \iota_2^{\varepsilon }\left( t \right) $ as follows,
\begin{eqnarray}\label{lemma3.37}
 \hat{\mathcal A}^\varepsilon \iota_2^{\varepsilon }\left( t \right) =  -  \sum\limits_{i,j = 1}^{n }{\iota_{{x_i}{x_j}}}\left( x \right) {\bar a^\nu_{ij}\left( {x,{y ^\varepsilon }\left( {t } \right)} \right) + \sum\limits_{i,j = 1}^{n }\iota_{{x_i}{x_j}}(x) \bar a_{ij}\left( {x} \right)}+O(\varepsilon).
\end{eqnarray}
Therefore, from (\ref{lemma3.18}) to (\ref{lemma3.37}),  define that
\begin{eqnarray*}\label{lemma3.511}
	\mathbb{E}\big\vert \int_{s}^{T} \big[{{\hat {\mathcal A}}^\varepsilon }{\iota^{\varepsilon }}\left( t \right)-{\mathcal A}\iota\left( {{x^{\varepsilon }}\left( t \right)} \right)\big]dt\big\vert 
	&\le&  {\mathcal{I}_1+\mathcal{I}_2+O(\varepsilon)},
\end{eqnarray*}
where
\begin{eqnarray*}
	\mathcal{I}_1
	&{\rm{ = }}& \mathbb{E}\Big\vert \int_{s}^{T}\big[  {{\iota_{{x}}}\left( x \right)f\left( {{x^{\varepsilon }}\left( t \right),{y ^\varepsilon  }\left( t \right),r^\varepsilon(t)} \right)} \cr
	&&\quad-{{\iota_{{x}}}\left( x \right)\bar f^\nu\left( {{x^{\varepsilon }}\left( t \right),{y ^\varepsilon }\left( {t } \right)} \right)} \big]dt\big\vert ,
\end{eqnarray*}
and 
\begin{eqnarray*}
	\mathcal{I}_2
	&{\rm{ = }}& \mathbb{E}\big\vert \frac{1}{2}\int_{s}^{T}\big[ {{\iota_{{x}{x}}}\left( x \right)a\left( {{x^{\varepsilon }}\left( t \right),{y ^\varepsilon  }\left( t \right),r^\varepsilon(t)} \right)} \cr
	&&\quad-  {{\iota_{{x}{x}}}\left( x \right)\bar a^\nu\left( {{x^{\varepsilon }}\left( t \right),{y ^\varepsilon }\left( {t} \right)} \right)}\big]dt \big\vert ,
\end{eqnarray*}
 {where $\iota_x(x)=\big(\frac{\partial {\iota(x)}}{\partial{x_1}},\frac{\partial {\iota(x)}}{\partial{x_2}},...,\frac{\partial {\iota(x)}}{\partial{x_n}}\big)$ and $\iota_{xx}(x)=\big[\frac{\partial^2 {\iota(x)}}{\partial{x_i}\partial{x_j}}\big]_{n\times n}$.} To be clear, we do not distinguish the dimension in the following parts.
Take  a partition of time interval $[s,T]$ of size $\delta$, when $k\in[{[ {{s \mathord{\left/{\vphantom {s \delta }} \right.\kern-\nulldelimiterspace} \delta }} ] }+1,{[ {{T \mathord{\left/{\vphantom {T \delta }} \right.\kern-\nulldelimiterspace} \delta }} ] }-1 ) $, 
\begin{eqnarray*}\label{lemma3.512}
	\mathcal{I}_1 &\le& \mathbb{E}\big\vert \sum\limits_{k = {\left[ {{s \mathord{\left/{\vphantom {s \delta }} \right.\kern-\nulldelimiterspace} \delta }} \right] }+1}^{\left[ {{T  \mathord{\left/{\vphantom {T \delta }} \right.	\kern-\nulldelimiterspace} \delta }} \right]-1}\int_{k\delta}^{\left( {k + 1} \right)\delta} \big[ {{\iota_{{x}}}\left( x \right)f\left( {{x^{\varepsilon }}\left( t \right),{y ^\varepsilon  }\left( t \right),r^\varepsilon(t)} \right)}-{{\iota_{{x}}}\left( x \right)\bar f^\nu\left( {{x^{\varepsilon }}\left( t \right),{y ^\varepsilon }\left( {t } \right)} \right)}\big]dt\big\vert \cr
	&&+\mathbb{E} {\Big\vert  {\int_{\left[ {{T \mathord{\left/{\vphantom {(T-u) \delta }} \right.	\kern-\nulldelimiterspace} \delta }} \right]\delta }^{{T} } {\big[ {{\iota_{{x}}}\left( x \right)f\left( {{x^{\varepsilon }}\left( t \right),{y ^\varepsilon  }\left( t \right),r^\varepsilon(t)} \right)}-{{\iota_{{x}}}\left( x \right)\bar f^\nu\left( {{x^{\varepsilon }}\left( t \right),{y ^\varepsilon }\left( {t } \right)} \right)} \big]dt} }  \Big\vert }\cr
	&&+\mathbb{E} {\Big\vert  {\int_{ s  }^{{\left(\left[ {{s \mathord{\left/{\vphantom {s \delta }} \right.\kern-\nulldelimiterspace} \delta }} \right] +1\right)\delta} } {\big[ {{\iota_{{x}}}\left( x \right)f\left( {{x^{\varepsilon }}\left( t \right),{y ^\varepsilon  }\left( t \right),r^\varepsilon(t)} \right)}-{{\iota_{{x}}}\left( x \right)\bar f^\nu\left( {{x^{\varepsilon }}\left( t \right),{y ^\varepsilon }\left( {t } \right)} \right)} \big]dt} }  \Big\vert }\cr
	&=:& \mathcal{I}_{11}+\mathcal{I}_{12}+\mathcal{I}_{13}.
\end{eqnarray*}
Then by the \lemref{lem3.4}, we have
\begin{eqnarray*}\label{lemma3.515}
	\mathcal{I}_{12} \le C \delta,\quad
	\mathcal{I}_{13} \le C \delta.
\end{eqnarray*}
Furthermore, compute the  term $ {\mathcal{I}_{11}}$,
\begin{eqnarray*}\label{lemma3.513}
	\mathcal{I}_{11} &\le& \sum\limits_{k = {\left[ {{s \mathord{\left/{\vphantom {s \delta }} \right.\kern-\nulldelimiterspace} \delta }} \right] }+1}^{\left[ {{T  \mathord{\left/{\vphantom {T \delta }} \right.	\kern-\nulldelimiterspace} \delta }} \right]-1}\mathbb{E}\big\vert \int_{k\delta}^{\left( {k + 1} \right)\delta} \big[ {{\iota_{{x}}}\left( x \right)f\left( {{x^{\varepsilon }}\left( t \right),{y ^\varepsilon  }\left( t \right),r^\varepsilon(t)} \right)}-{{\iota_{{x}}}\left( x \right)\bar f^\nu\left( {{x^{\varepsilon }}\left( t \right),{y ^\varepsilon }\left( {t } \right)} \right)}\big]dt\big\vert \cr
	&\le&\left(\left[ {{T \mathord{\left/{\vphantom {t \delta }} \right.\kern-\nulldelimiterspace} \delta }} \right]-\left[ {{s \mathord{\left/{\vphantom {s \delta }} \right.\kern-\nulldelimiterspace} \delta }} \right]\right)\mathop {\max }\limits_{{{\left[ {{s \mathord{\left/{\vphantom {u \delta }} \right.\kern-\nulldelimiterspace} \delta }} \right] }+1} \le k \le \left[ {{T \mathord{\left/{\vphantom { \delta }} \right.\kern-\nulldelimiterspace} \delta }} \right]-1}\mathbb{E} \Big\vert  \int_{0}^{{\delta} } \big[ {{\iota_{{x}}}\left( x \right)f\left( {{x^{\varepsilon }}\left( t \right),{y ^\varepsilon  }\left( t \right),r^\varepsilon(t)} \right)}\big.\Big.\cr
	&&\Big.\big.-{\iota_{{x}}}\left( x\left( k \delta \right) \right){f\left( {{x^{\varepsilon }}\left( k \delta \right),{y ^\varepsilon  }\left( t \right),r^\varepsilon(t)} \right)} \big]dt   \Big\vert \cr
	&&+\left(\left[ {{T \mathord{\left/{\vphantom {t \delta }} \right.\kern-\nulldelimiterspace} \delta }} \right]-\left[ {{s \mathord{\left/{\vphantom {u \delta }} \right.\kern-\nulldelimiterspace} \delta }} \right]\right)\mathop {\max }\limits_{{{\left[ {{s \mathord{\left/{\vphantom {u \delta }} \right.\kern-\nulldelimiterspace} \delta }} \right] }+1} \le k \le \left[ {{T \mathord{\left/{\vphantom { \delta }} \right.\kern-\nulldelimiterspace} \delta }} \right]-1}\mathbb{E} \Big\vert  \int_{0}^{{\delta} } {\iota_{{x}}}\left( x(k \delta) \right)\big[ {f\left( {{x^{\varepsilon }}\left( k \delta \right),{y ^\varepsilon  }\left( t \right),r^\varepsilon(t)} \right)}\big.\Big.\cr
	&&\Big.\big.-{\bar f^\nu\left( {{x^{\varepsilon }}\left( k \delta \right),{y ^\varepsilon }\left( {t } \right)} \right)} \big]dt   \Big\vert \cr
	&&+\left(\left[ {{T \mathord{\left/{\vphantom {t \delta }} \right.\kern-\nulldelimiterspace} \delta }} \right]-\left[ {{s \mathord{\left/{\vphantom {s \delta }} \right.\kern-\nulldelimiterspace} \delta }} \right]\right)\mathop {\max }\limits_{{{\left[ {{s \mathord{\left/{\vphantom {u \delta }} \right.\kern-\nulldelimiterspace} \delta }} \right] }+1} \le k \le \left[ {{T \mathord{\left/{\vphantom { \delta }} \right.\kern-\nulldelimiterspace} \delta }} \right]-1}\mathbb{E} \Big\vert  \int_{0}^{{\delta} } \big[ {{\iota_{{x}}}\left( x(k \delta) \right)\bar f^\nu\left( {{x^{\varepsilon }}\left( k \delta \right),{y ^\varepsilon }\left( {t } \right)} \right)}\big.\Big.\cr
	&&\Big.\big.-{{\iota_{{x}}}\left( x(t) \right)\bar f^\nu\left( {{x^{\varepsilon }}\left( t \right),{y ^\varepsilon }\left( {t } \right)} \right)} \big]dt   \Big\vert \cr
	&\le&\left(\left[ {{T \mathord{\left/{\vphantom {t \delta }} \right.\kern-\nulldelimiterspace} \delta }} \right]-\left[ {{s \mathord{\left/{\vphantom {s \delta }} \right.\kern-\nulldelimiterspace} \delta }} \right]\right)\mathop {\max }\limits_{{{\left[ {{s \mathord{\left/{\vphantom {u \delta }} \right.\kern-\nulldelimiterspace} \delta }} \right] }+1} \le k \le \left[ {{T \mathord{\left/{\vphantom { \delta }} \right.\kern-\nulldelimiterspace} \delta }} \right]-1}\mathbb{E} \Big\vert  \int_{0}^{{\delta} } {\iota_{{x}}}\left( x (k \delta)\right)\big[ {f\left( {{x^{\varepsilon }}\left( k \delta \right),{y ^\varepsilon  }\left( t \right),r^\varepsilon(t)} \right)}\big.\Big.\cr
	&&\Big.\big.-{\bar f^\nu\left( {{x^{\varepsilon }}\left( k \delta \right),{y ^\varepsilon }\left( {t } \right)} \right)} \big]dt   \Big\vert +C O(\delta)\cr
	&=:&\left(\left[ {{T \mathord{\left/{\vphantom {t \delta }} \right.\kern-\nulldelimiterspace} \delta }} \right]-\left[ {{s \mathord{\left/{\vphantom {s \delta }} \right.\kern-\nulldelimiterspace} \delta }} \right]\right)\mathop {\max }\limits_{{{\left[ {{s \mathord{\left/{\vphantom {u \delta }} \right.\kern-\nulldelimiterspace} \delta }} \right] }+1} \le k \le \left[ {{T\mathord{\left/{\vphantom { \delta }} \right.\kern-\nulldelimiterspace} \delta }} \right]-1}\mathcal{I}^{11}_1+C O(\delta),
\end{eqnarray*}
 {where the third inequality is from the (\ref{lemma3.0}), $\iota \left(  \cdot  \right)\in C_0^4\left( {{\mathbb{R}^n},\mathbb{R}} \right)$, and Lipschitz continuity of function $f$.} With aid of the \lemref{lem3.1} and the fact of independence between the Markovian switching and fast component (\ref{lemma3.11}), it   {obtains} 
\begin{eqnarray*}\label{lemma3.514}
	\mathcal{I}^{11}_1 &\le& \Big\{\mathbb{E} \Big\vert  {\int_{0}^{{\delta} } {{\iota_{{x}}}\left( x (k \delta)\right)\big[ {f\left( {{x^{\varepsilon }}\left( k \delta \right),{y ^\varepsilon  }\left( t \right),r^\varepsilon(t)} \right)}-{\bar f^\nu\left( {{x^{\varepsilon }}\left( k \delta \right),{y ^\varepsilon }\left( {t } \right)} \right)} \big]dt} }  \Big\vert ^2\Big\}^{{1 \mathord{\left/{\vphantom {1 2}} \right.\kern-\nulldelimiterspace} 2}}\cr
	&\le& \Bigg\{ \int_0^\delta \int_\tau ^\delta {{\Big\{ {{\mathbb{E}^{y,r_0}}{{\left\{ {\left[ {f\left( {{x^\varepsilon }\left( {k\delta } \right),{{  {y^\varepsilon}}}\left( \tau  \right)},{r^{\varepsilon }}\left( {\tau  } \right) \right) - \bar f^\nu\left( {{x^\varepsilon }\left( {k\delta } \right)},{{  {y^\varepsilon}}}\left( \tau  \right) \right)} \right]} \right\}}^2}} \Big\}}^{{1 \mathord{\left/{\vphantom {1 2}} \right.\kern-\nulldelimiterspace} 2}}}\Bigg.\cr
	&&\quad\cdot\Bigg.\Big\{ {\mathbb{E}^{y,r_0}}\left\{ {\mathbb{E}^{{y\left( \tau  \right),r(\tau)}}}\left\{ \left[ f\left( {{x^\varepsilon }\left( {k\delta } \right),{{ y^\varepsilon}}\left( {t - \tau } \right)},{r^{\varepsilon }}\left(t- {\tau  } \right) \right)\right.\right.\right.\Big.\Bigg.\cr
	&&\qquad\Bigg.\Big.\left.\left.\left.- \bar f^\nu\left( {{x^\varepsilon }\left( {k\delta } \right)} ,{{ y^\varepsilon}}\left( {t - \tau } \right)\right) \right] \right\} \right\}^2 \Big\}^{{1 \mathord{\left/{\vphantom {1 2}} \right.\kern-\nulldelimiterspace} 2}}  dtd\tau \Bigg\}^{ {1/2}}\cr
	&\le& \Bigg\{ \int_0^\delta \int_\tau ^\delta {{\Big\{ {{\mathbb{E}^{y,r_0}}{{\left\{ {\left[ {f\left( {{x^\varepsilon }\left( {k\delta } \right),{{ y^\varepsilon}}\left( \tau  \right)},{r^{\varepsilon }}\left( {\tau  } \right) \right) - \bar f^\nu\left( {{x^\varepsilon }\left( {k\delta } \right)},{{ y^\varepsilon}}\left( \tau  \right) \right)} \right]} \right\}}^2}} \Big\}}^{{1 \mathord{\left/{\vphantom {1 2}} \right.\kern-\nulldelimiterspace} 2}}}\Bigg.\cr
	&&\quad\times\Bigg.\Big\{ {\mathbb{E}^{y,r_0}}\Big\{ {\sum\limits_{\gamma = 1}^n }\int_{{\mathbb{R}^n}}\left\{  f\left( {x^\varepsilon }\left( {k\delta } \right),{{ \xi}},\gamma \right){P}\left( {{y ^\varepsilon }\left( \tau \right),\tau,d\xi,t} \right)\big[{P}\left( r^\varepsilon(t)=\gamma \vert r^\varepsilon(\tau)={\gamma}'  \right)\big.\right.\Big.\Big.\Bigg.\cr
	&&\qquad\Bigg.\Big.\Big.\left.\big.-{\nu_\gamma }\big]   \right\} \Big\}^2 \Big\}^{{1 \mathord{\left/{\vphantom {1 2}} \right.\kern-\nulldelimiterspace} 2}}dtd\tau\Bigg\}^{ {1/2}} \cr
	&\le&C \varepsilon.
\end{eqnarray*}
Similarly, it can prove that when $\varepsilon$ tends to zero, the term $\mathcal{I}_2$ tends to zero.

 {Step 3: With the tightness of $\left\{ {{x^\varepsilon }\left(  \cdot  \right)} \right\}$ in \lemref{lem3.4}, it gets that {${x^\varepsilon }\left(  \cdot  \right)$ converges weakly to  $x\left(  \cdot  \right)$}.}

\begin{rem}\label{rem5.10}
	For slow-fast system without Markovian switching regimes as follows, 
	\begin{eqnarray*}
		\begin{cases}
			d{x^\varepsilon}\left(t\right)&=f\left({{x^\varepsilon}\left(t\right),{y^\varepsilon}\left(t\right)}\right)dt+\sigma\left({{x^\varepsilon }\left(t\right),{y^\varepsilon}\left(t\right)} \right)dW_1\left(t\right)\cr
			&\quad +\int_{ \vert z \vert <c}{g\left({{x^\varepsilon}\left(t\right),z}\right)\tilde N_1\left( {dt,dz} \right)},\cr
			d{y^\varepsilon}\left(t\right)&=\frac{1}{\varepsilon }{b}\left( {{y^\varepsilon}\left(t\right)}\right)dt + \frac{1}{{\sqrt \varepsilon  }}V \left({{y ^\varepsilon }\left( t \right)} \right)d{W_2}\left(t\right)+\int_{ \vert z \vert <c}{h\left({{y ^\varepsilon}\left(t\right),z}\right)\tilde N_2^\varepsilon \left({dt,dz}\right)},
		\end{cases}
	\end{eqnarray*} 
	under the Assumptions {\rm (A1)} to {\rm (A5)}, \thmref{thm3.4} can be  {obtained}.
\end{rem}

\section{Examples}\label{Sec-4}
\begin{exm}\label{exm1}
	 {Consider the  stochastic vibration with damping, when  the drag is big enough, }
	\begin{eqnarray}\label{sec4-1}
	\begin{cases}
	d{x^\varepsilon}\left(t\right)&={y^\varepsilon}\left(t\right)dt,\cr
	d{y^\varepsilon}\left(t\right)&=-\frac{1}{\varepsilon } {{y^\varepsilon}\left(t\right)}dt+\frac{1}{\varepsilon } dt + \frac{1}{{\sqrt \varepsilon  }}  d{W}\left(t\right),
	\end{cases}
	\end{eqnarray}
	where $x^\varepsilon(0)=x_0, y^\varepsilon(0)=y_0$.  {The small parameter $\varepsilon$ is related to the     drag. The slow component and fast component are motion and speed respectively. }
\end{exm}
 Therefore, $\mu$ is the unique invariant measure of the transition semigroup for the  following system,
\begin{eqnarray}\label{sec4-10}
	dy=-y dt  + dt+ dW\left( t \right).
\end{eqnarray}
According to the Fokker-Planck-Kolmogorov equation \cite{Xu2019Path}, the exact solution of its transition density function is as follows,
\begin{eqnarray*}
	p(y, t\vert y_{0},0)=\sqrt{\frac{1}{2 \pi D\left(1-e^{-2 t}\right)}} \exp \left[-\frac{[y-1-(y_{0})e^{-t}]^{2}}{2 D\left(1-e^{-2 t}\right)}\right],
\end{eqnarray*}
where $D$ represents the noise strength, and the stationary solution is
\begin{eqnarray*}
	\mu(y)=\sqrt{\frac{1}{2 \pi D}} \exp \left\{-\frac{(y-1) ^{2}}{2 D}\right\}.
\end{eqnarray*} 
Then the averaged system is as follows,
\begin{eqnarray}\label{sec4-11}
{\bar x}\left(t\right)&=x_0+t,
\end{eqnarray}
the simulation for expectation of  {solutions} to original slow system (\ref{sec4-1}) and averaged system (\ref{sec4-11}) is in Figure 1.
\begin{figure}[H]
	\centering
	{\subfigure[]{
			\includegraphics[width=0.35\textwidth,height=0.3\textwidth]{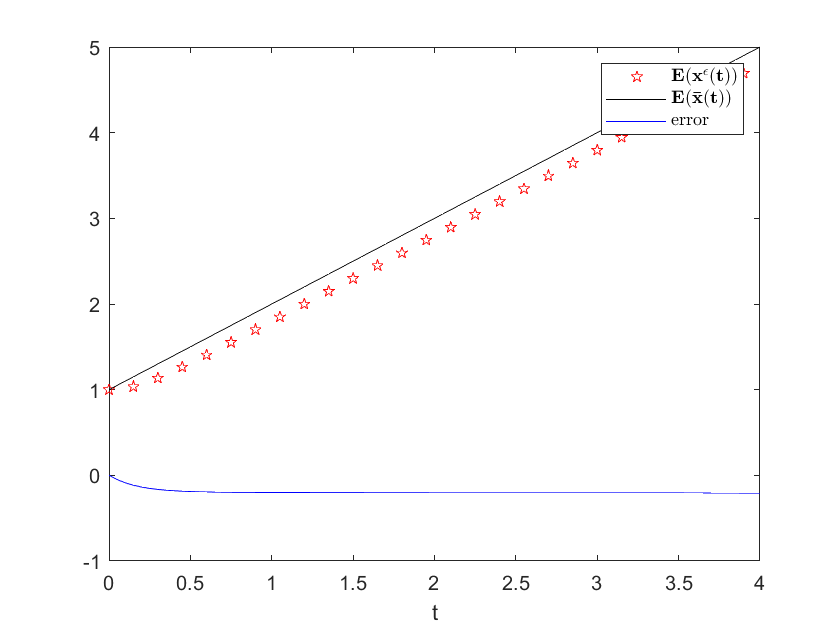}}}
	\hspace{0.2\textwidth}
	\subfigure[]{
		\includegraphics[width=0.35\textwidth,height=0.3\textwidth]{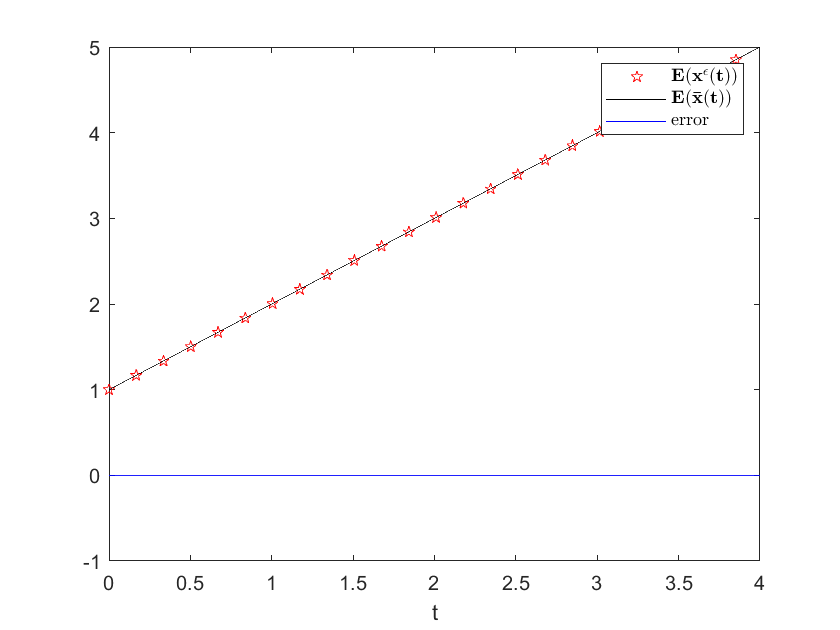}}
	\caption{ {Expectation of solutions $x^\varepsilon(t)$, $\bar x(t)$ to the original slow system (\ref{sec4-1}) and the averaged system  (\ref{sec4-11}), and error curve. $(a)$ $D=0.1$, $x_0=1,y_0=0$,  $\varepsilon=0.01$,  $(b)$ $D=0.1$, $x_0=1,y_0=0$, $\varepsilon=0.001$.} }
	\label{fig1}
\end{figure}
  {In Fig.1(a), error curves  show that the expectation of solutions to  averaged equation (\ref{sec4-11}) has a good approximation of those to the original slow equation  with  $\varepsilon=0.01$. In Fig.1(b),   similarly  for the expectation of solutions to the original slow system (\ref{sec4-1}) and the averaged system (\ref{sec4-11})   with  $\varepsilon=0.001$,  the error becomes smaller with smaller $\varepsilon$.  Then, the main \thmref{thm3.4} can be verified in this stochastic vibration model with damping.}

\begin{exm}\label{exm2}
	\begin{eqnarray}\label{sec4-2}
	\begin{cases}
	d{x^\varepsilon}\left(t\right)&=-f\left(r^\varepsilon(t),{x^\varepsilon}\left(t\right),{y^\varepsilon}\left(t\right)\right)dt-\sigma(r^\varepsilon(t),{y^\varepsilon}\left(t\right))dW_1\left(t\right),\cr
	d{y^\varepsilon}\left(t\right)&=-\frac{1}{\varepsilon } {{y^\varepsilon}\left(t\right)}dt + \frac{1}{{\sqrt \varepsilon  }}  d{W_2}\left(t\right),
	\end{cases}
	\end{eqnarray}
	where  $x^\varepsilon(0)=x_0, y^\varepsilon(0)=y_0$, $f(1,x^\varepsilon(t),{y^\varepsilon}\left(t\right))=-(x^\varepsilon(t)+{y^\varepsilon}\left(t\right))$ for $r^\varepsilon(t)=1$ and $f(2,x^\varepsilon(t),{y^\varepsilon}\left(t\right))=-2(x^\varepsilon(t)+{y^\varepsilon}\left(t\right))$ for $r^\varepsilon(t)=2$, and $\sigma(1,{y^\varepsilon}\left(t\right))=-{y^\varepsilon}\left(t\right)$ for $r^\varepsilon(t)=1$ and $\sigma(2,{y^\varepsilon}\left(t\right))=-2{y^\varepsilon}\left(t\right)$ for $r^\varepsilon(t)=2$.
\end{exm}
Therefore, $\mu$ is the unique invariant measure of the transition semigroup for the  system as follows:
\begin{eqnarray*}
	dy=-y dt  + dW_2\left( t \right).
\end{eqnarray*}
According to the Fokker-Planck-Kolmogorov equation, 
\begin{eqnarray*}
	\mu(y)=\sqrt{\frac{1}{2 \pi D}} \exp \left\{-\frac{y ^{2}}{2 D}\right\}.
\end{eqnarray*}
Then, $r^\varepsilon(t)$ is a fast varying Markov chains whose generator is given by
\begin{eqnarray*}
	{{\mathcal{Q}}^\varepsilon }: = \frac{{\widetilde {\mathcal{Q}}}}{\varepsilon },
\end{eqnarray*}
where 
\begin{equation}       
{\widetilde {\mathcal{Q}}}=\left(                 
\begin{array}{ccc}   
-1 & 1 \\  
2 & -2 \\  
\end{array}
\right)                 
\end{equation}\nonumber
has a unique solution $\nu = \left(\frac{1}{2 },\frac{1}{2 } \right)$ which is termed a quasi-stationary distribution. 
Then the  averaged system is as follows,
\begin{eqnarray}\label{sec4-12}
d{\bar x}\left(t\right)&=-\frac{3}{2 }{\bar x}\left(t\right)dt-\sqrt{\frac{5}{2 }}dW_1\left(t\right),
\end{eqnarray}
then the simulation for expectation of  {solutions} to the original slow system (\ref{sec4-2}) and the averaged system (\ref{sec4-12}) is in Figure 2.
\begin{figure}[H]
	\centering
	{\subfigure[]{
			\includegraphics[width=0.35\textwidth,height=0.3\textwidth]{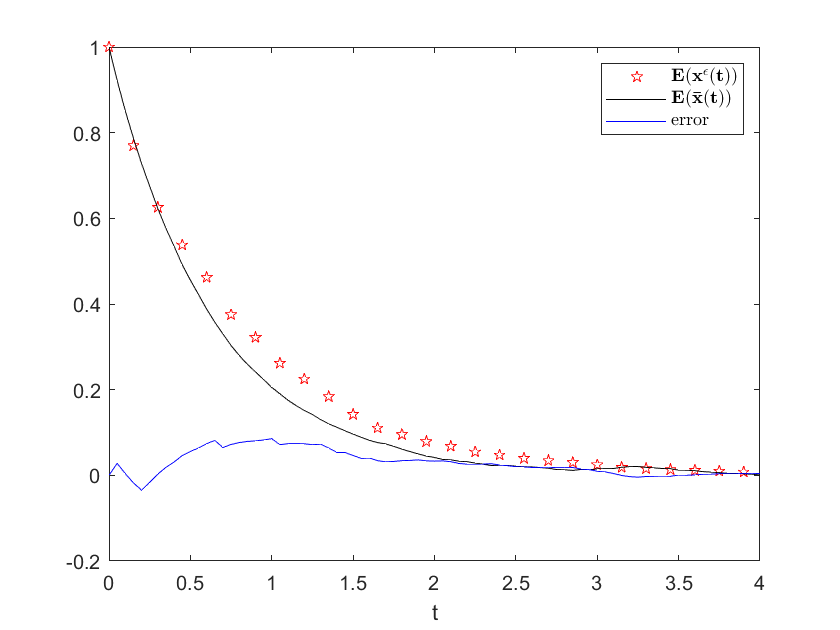}}}
	\hspace{0.2\textwidth}
	\subfigure[]{
		\includegraphics[width=0.35\textwidth,height=0.3\textwidth]{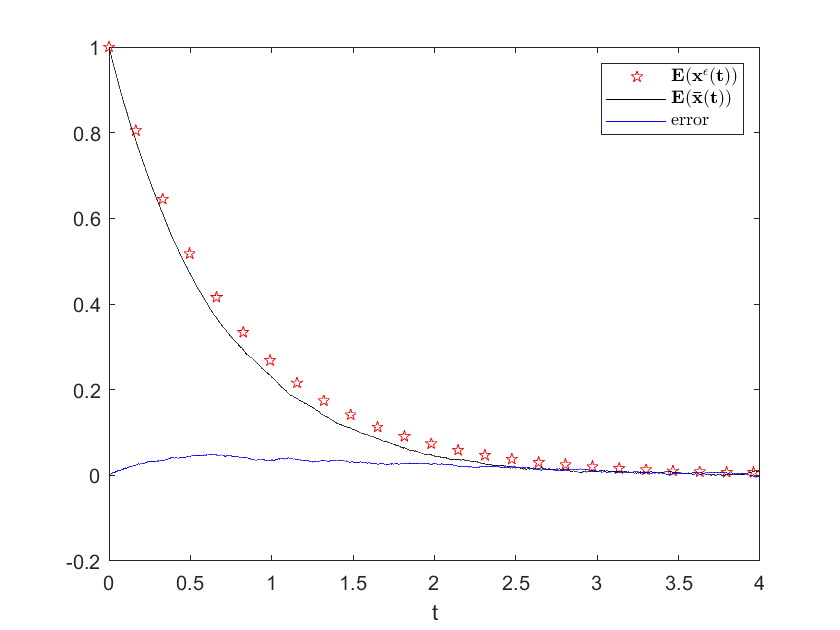}}
	\caption{ {Expectation of solutions $x^\varepsilon(t)$, $\bar x(t)$ to the original slow system (\ref{sec4-2}) and the averaged system  (\ref{sec4-12}), and error curve. $(a)$ $D=0.1$, $x_0=1,y_0=0$,  $\varepsilon=0.2$,  $(b)$ $D=0.1$, $x_0=1,y_0=0$, $\varepsilon=0.001$.} }
	\label{fig2}
\end{figure}
  {Example 4.2 considers the fast component and switching regimes simultaneously, in which the exact solution  is difficult to be obtained. The averaged equation (\ref{sec4-12}) is a 1-dimensional linear SDE, while it is simpler with its solution to be easily solved  \cite{Arnold1974Stochastic}. In Fig.2(a) and (b),  a good agreement between solutions to the averaged system (\ref{sec4-12}) and the original system (\ref{sec4-2})  is observed.  It implied that the averaged solution (\ref{sec4-12}) can approximate the original solution to the slow-fast system (\ref{sec4-2}), and the weak convergence is assured by \thmref{thm3.4}. }

\begin{exm}\label{exm3}
	\begin{eqnarray}\label{sec4-3}
	\begin{cases}
	d{x^\varepsilon}\left(t\right)&=-f\left(r^\varepsilon(t),{x^\varepsilon}\left(t\right),{y^\varepsilon}\left(t\right)\right)dt-\sigma(r^\varepsilon(t),{y^\varepsilon}\left(t\right))dW_1\left(t\right)\cr
	&\quad+\int_{ \vert z \vert <c}{{{x^\varepsilon}\left(t\right)z}\tilde N_1\left( {dt,dz} \right)},\cr
	d{y^\varepsilon}\left(t\right)&=-\frac{1}{\varepsilon } {{y^\varepsilon}\left(t\right)}dt + \frac{1}{{\sqrt \varepsilon  }}  d{W_2}\left(t\right),
	\end{cases}
	\end{eqnarray}
	where  $x^\varepsilon(0)=x_0, y^\varepsilon(0)=y_0$,  $f(r^\varepsilon(t),x^\varepsilon(t),{y^\varepsilon}\left(t\right))$ and $\sigma(r^\varepsilon(t),{y^\varepsilon}\left(t\right))$ are same  in Example 4.2. 
\end{exm}
There exist jumps due to the compound Poisson random measure, 
where $dJ(t)=\int_{ \vert z \vert <c}{{z}\tilde N\left( {dt,dz} \right)}$. $J(t):=\sum_{j=1}^{N(t)}Z_j$ where $Z_j$ are independent identically distributed normal random variables, $N(t)$ denote the Poisson process with intensity $\lambda>0$.

With the stationary solution of  its transition density function for the system (\ref{sec4-10}),
then the  averaged system is 
\begin{eqnarray}\label{sec4-31}
d{\bar x}\left(t\right)&=-\frac{3}{2 }{\bar x}\left(t\right)dt-\sqrt{\frac{5}{2 }}dW_1\left(t\right)+\int_{ \vert z \vert <c}{{{\bar x}\left(t\right)z}\tilde N_1\left( {dt,dz} \right)},
\end{eqnarray}
 the simulation for expectation of  {solutions} to the original slow system (\ref{sec4-3}) and the averaged system (\ref{sec4-31}) is in Figure 3.
\begin{figure}[H]
	\centering
	{\subfigure[]{
			\includegraphics[width=0.35\textwidth,height=0.3\textwidth]{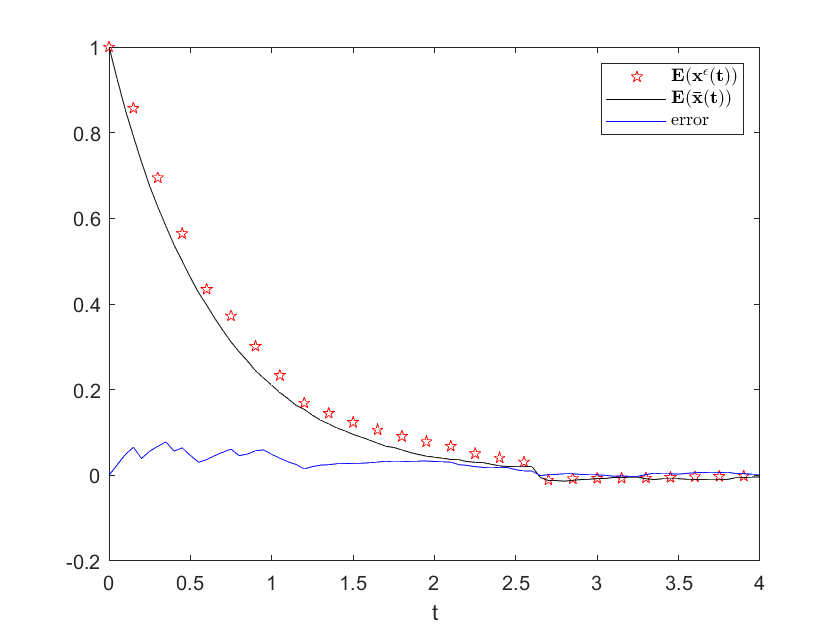}}}
	\hspace{0.2\textwidth}
	\subfigure[]{
		\includegraphics[width=0.35\textwidth,height=0.3\textwidth]{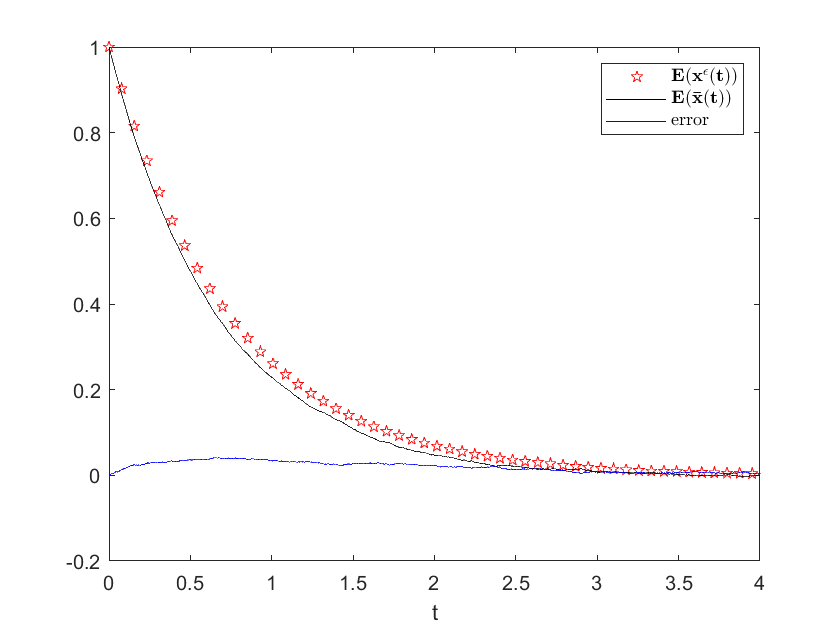}}
	\caption{ {Expectation of solutions $x^\varepsilon(t)$ and $\bar x(t)$ to the original slow system (\ref{sec4-3}) and the averaged system (\ref{sec4-31}), and error curve. $(a)$ $D=0.1$,  $\lambda=0.5$,  and $Z_j \sim N(0,1)$, $x_0=1,y_0=0$, $\varepsilon=0.2$, $(b)$ $D=0.1$,  $\lambda=0.5$,  and $Z_j \sim N(0,1)$, $x_0=1,y_0=0$, $\varepsilon=0.001$. }}
	\label{fig3}
\end{figure}
  {Based on the previous analysis, Example 4.3 considers the hybrid slow-fast system with jumps. In Fig.3(a) and (b), one can see that, solutions to   the averaged system (\ref{sec4-31}) still have a good agreement with those to the original system (\ref{sec4-3}) under different small parameters $\varepsilon$, that coincides with the theoretical results in \thmref{thm3.4}.}

\section{Concluding Remarks}\label{Sec-5}

 In this paper,  weak convergence for  slow-fast systems with jumps modulated by Markovian switching  is investigated by martingale method.  {To eliminate the fast component and Markovian switching in averaging, the new approach, a combination of perturbed test functions and time discretization, is applied efficiently.} This causes a computation burden in our proof. Therefore, other methods which could avoid this need to be developed in the future. Then we point out that the fast component is independent of the slow component. Since the coupling between the slow and fast components is widespread, the case that the fast component is dependent on the slow component is worthy to be noticed. However, if the fast component depends on the slow component,  we have to consider the transition probability where the slow component is as a parameter. This is an obstacle when applying this martingale method, which will be considered in  future work.

\section*{Acknowledgements}
Y. Xu was partly supported by the Key International
(Regional) Cooperative Research Projects of the NSF of China (Grant
12120101002), the NSF of China (Grant 12072264), the Fundamental
Research Funds for the Central Universities, the Research Funds
for Interdisciplinary Subject of Northwestern Polytechnical University, the
Shaanxi Provincial Key R\&D Program (Grants 2020KW-013, 2019TD-010). X. Yang was partly supported by the NSF of China (Grant 11802236).
B. Pei was partially supported by the NSF of China (Grant 11802216), the NSF of
Chongqing (Grant cstc2021jcyj-msxmX0296), the Guangdong Basic and
Applied Basic Research Foundation, the Fundamental Research Funds for the Central
Universities, the Young Talent fund of University Association for Science and
Technology in Shaanxi, China.
\section*{Availability of data and material}
Not applicable.

\section*{Competing interests}
The authors declare that they have no competing interests.

\section*{Authors' contributions}
All authors completed the paper together. All authors read and approved the final manuscript.

\end{document}